\documentstyle[amscd]{amsart}
\makeatletter
\renewcommand{\subsection}{\@startsection{subsection}{2}{\z@}%
{\baselineskip}{0.5\baselineskip}{\defaultfont\bf}}
\makeatother
\newtheorem{lem}{Lemma}[section]
\newtheorem{pro}[lem]{Proposition}

\theoremstyle{definition}
\newtheorem{defn}{Definition}[section]
\numberwithin{equation}{section}
\def\dj{d\kern-.30em\raise1.25ex\vbox{\hrule width .3em height .03em}}
\def\Dj{D\kern-.75em\raise0.75ex\vbox{\hrule width .3em height .03em}
\kern.03em}
\newsymbol\restr 1316
\newsymbol\varkappa 207B
\newcommand{\ADP}{\cal{C}(P)}
\newcommand{\AD}{\cal{C}}
\newcommand{\SC}{S_c}
\newcommand{\WC}{\Omega_c}
\newcommand{\hor}{\mbox{\family{euf}\shape{n}\selectfont hor}}
\newcommand{\grten}{\mathbin{\widehat{\otimes}}}
\newcommand{\ad}{\mbox{\shape{n}\selectfont ad}}
\newcommand{\id}{\mbox{\shape{n}\selectfont id}}
\newcommand{\Ad}{\varpi}
\newcommand{\lie}{\mbox{\family{euf}\shape{n}\selectfont lie}}
\newcommand{\e}{\epsilon}
\newcommand{\k}{\kappa}
\newcommand{\Sum}{\displaystyle{\sum}}
\newcommand{\ver}{\mbox{\family{euf}\shape{n}\selectfont ver}}
\newcommand{\con}{\mbox{\family{euf}\shape{n}\selectfont con}}
\newcommand{\inv}{i\!\hspace{0.8pt}n\!\hspace{0.6pt}v}
\begin{document}
\title[Quantum Principal Bundles]{Quantum Principal Bundles and\\
Corresponding Gauge Theories}
\author{Mi\'co \Dj ur\Dj evi\'c}
\address{Instituto de Matematicas, UNAM, Area
de la Investigacion Cientifica, Circuito Exterior,
Ciudad Universitaria, M\'exico DF, CP 04510, MEXICO\newline
\indent {\it Written In}\newline
\indent Faculty of Physics, University of Belgrade, Belgrade, SERBIA\newline
\indent{\it Final Version}\newline
\indent Facultad de Estudios Superiores Cuautitlan, UNAM, MEXICO}
\maketitle
\begin{abstract} A  generalization  of  classical  gauge  theory
is presented,  in  the  framework   of   a   noncommutative-geometric
formalism of  quantum  principal  bundles  over  smooth  manifolds.
Quantum  counterparts  of
classical gauge  bundles,  and  classical  gauge  transformations,
are introduced and investigated. A natural  differential  calculus
on quantum gauge bundles is constructed and analyzed.  Kinematical
and dynamical  properties  of  corresponding  gauge  theories  are
discussed.
\end{abstract}
\tableofcontents
\section{Introduction}
\renewcommand{\thepage}{}
     The aim of this study  is  to  present  a  generalization  of
classical gauge theory, in which quantum groups play the  role  of
entities describing local symmetries.

     All  considerations  will  be  performed  within  a   general
conceptual framework of non-commutative differential geometry \cite{C}.

     The whole paper is based on a noncommutative-geometric theory
of principal bundles over classical smooth  manifolds,  possessing
quantum structure groups. This theory is presented in
\cite{D}. Here, fundamental  structural  elements  of  classical  gauge
theory will be generalized and incorporated into the formalism  of
quantum principal bundles.

     The paper is organized as follows.

     In the next section, a preparatory material is collected.  As
first, we fix the notation and  introduce  in  the  game  relevant
quantum group entities. Secondly, we present  the  most  important
ideas and  results  of  \cite{D},  which  will  be  used  in  the  main
considerations.

     The starting point for all constructions of this paper  is  a
quantum principal $G$-bundle $P$ over  a  smooth   manifold   $M$.
Here,  $M$
plays the role of space-time  while  $G$  is  a  compact  matrix
quantum (structure) group \cite{W2} representing ``local  symmetries''  of
the system.

     In Section 3  a quantum analogue of the gauge bundle will be
constructed and investigated. This quantum bundle (over $M$) will be
denoted  by  $\ADP$.  Various   quantum   counterparts   of   gauge
transformations  are  naturally  associated  to  $\ADP$. Further, a
differential calculus on the bundle $\ADP$ will be constructed,  by
combining the  standard  differential  calculus  on  $M$  (based  on
differential forms) with an appropriate differential  calculus  on
the quantum group  $G$.  This  calculus  on  $\ADP$  is  relevant  in
situations in which quantum counterparts of gauge  transformations
act on entities related to differential calculus on the  principal
bundle $P$ (connection forms, for example).

     It is important to  mention  that  there  exist  two  natural
inequivalent ways of introducing  quantum  counterparts  of  gauge
transformations. The first one is to translate  into  the  quantum
context  the  idea  that  gauge   transformations   are   vertical
automorphisms of the principal bundle $P$. This approach leads to  a
standard group (of gauge transformations of  $P$).  The  same  group
will be obtained if we consider counterparts of  sections  of  the
bundle $\ADP$. However, it turns out that such a concept of a gauge
transformation does not describe gauge-like phenomenas related  to
the quantum nature of the space $G$. Namely, because of the inherent
geometrical  inhomogeneity  of  quantum   groups,   every   quantum
principal  bundle  $P$  over  $M$  is  completely  determined  by  its
classical part $P_{cl}$   (interpretable as the set of points of $P$).  The
classical part is an ordinary principal $G_{cl}$-bundle over  $M$,  where
$G_{cl}$  is  a  group  (the  classical  part  of  $G$)  interpretable  as
consisting  of  points  of  $G$.   We   shall   prove   that   gauge
transformations of $P$ are in  a  natural  bijection  with  standard
gauge transformations of $P_{cl}  $. Further, we shall prove that the set
of points of $\ADP$  coincides,  in  a  natural  manner,  with  the
standard gauge bundle $\AD(P_{cl})$.
\renewcommand{\thepage}{\arabic{page}}
The second approach to gauge transformations  is  in  some  sense
indirect. The main idea is to construct the ``action'' of the bundle
$\ADP$ on $P$ (generalizing the classical situation).  This  approach
does not meet geometrical obstacles. In  classical  geometry,  the
mentioned action naturally  contains  all  the  information  about
gauge transformations.

     Section 4 is devoted to the formulation and kinematical  and
dynamical  analysis  of  quantum  group  gauge  theories,  in  the
framework of quantum  principal  bundles.  Gauge  fields  will  be
geometrically represented by connections on $P$. Internal degrees of
freedom  of  such  gauge  fields  are  determined  by   fixing   a
bicovariant  first-order  differential *-calculus \cite{W3}   on  the
structure quantum group $G$. In this paper  we  shall  deal  with  a
unique differential calculus on $G$ which can  be  characterized  as
the  minimal  bicovariant  differential  calculus  compatible,  in
appropriate sense, with the geometrical structure  on  the  bundle
$P$.
If we start from this calculus on the group then it is possible to
built natural differential calculi on bundles $P$  and  $\ADP$  which
are always ``locally trivialized'' when  bundles  $P$ (and  $\ADP$) are
locally trivialized.

     Dynamical properties of the gauge theory will  be  determined
after fixing  an  appropriate  lagrangian.  In  analogy  with  the
classical gauge theory, we shall consider  lagrangians  which  are
quadratic functions of the curvature form. We  shall  compute  the
corresponding equations of motion. Symmetry properties  of  the
introduced  lagrangian
will be analyzed. We shall prove the invariance of the  lagrangian
under the action of the (ordinary) group of gauge  transformations
of $P$. Further, it turns out that the lagrangian is  invariant,  in
an appropriate sense, under the natural action of $\ADP$ on $P$. This
corresponds to the full gauge  invariance  of  the  lagrangian  in  the
classical theory.

 In Section 5 everything will  be  illustrated  on  a  simple  but
highly non-trivial example in which $G$ is the quantum $SU(2)$  group.
The most important observation is  that  the  corresponding  gauge
theory is {\it essentially different} from  the  classical
$SU(2)$  gauge
theory, and does not reduce  to  the  classical  theory  when  the
deformation parametar $1-\mu$ tends to zero. This is caused by the  fact
that  the  minimal  admissible  bicovariant  calculus  does   not
respect the classical limit.  Namely,  a  detailed
analysis \cite{D} shows that for $\mu\in (-1,1)\setminus\{0\}$
the space of left-invariant
elements (playing the role of the dual space of the  corresponding
Lie algebra) of  the  mentioned  minimal  calculus  is  infinitely
dimensional, and can be naturally identified with the  algebra  of
polynomial  functions  on   a   quantum   2-sphere.   Hence,   the
corresponding  gauge  fields  possess  infinitely  many   internal
degrees of freedom, in contrast to the classical case.
Finally, in Section 6 concluding remarks are made.

     The paper ends with an  Appendix,  in  which  some  technical
properties related to the minimal admissible bicovariant  calculus
on the quantum $SU(2)$ group are collected.

\section{ Mathematical Background}
     Let $G$ be a compact matrix quantum group \cite{W2}. We shall denote
by $\cal{A}$  the  *-algebra  of  ``polynomial  functions''  on  $G$,  and  by
$\phi\colon \cal{A}\rightarrow \cal{A}\otimes \cal{A},$
$\e\colon \cal{A}\rightarrow \Bbb{ C}$  and  $\k\colon  \cal{A}\rightarrow
\cal{A}$ the
comultiplication, counit  and the  antipode
respectively. The symbols $a^{(1)} \otimes \dots \otimes a^{(n)}$
will be used
for the
result  of  an   $(n-1)$-fold   comultiplication   of   an   element
$a\in \cal{A}$  (so
that $\phi (a)=a^{(1)} \otimes a^{(2)} )$. Let $G_{cl}$   be the
classical  part \cite{D} of  $G$.
Explicitly, $G_{cl}$ is  consisting  of  *-characters
(nontrivial multiplicative
linear hermitian functionals) of $\cal{A}$. The Hopf algebra structure  on
$\cal{A}$ naturally induces the group structure on $G_{cl},$ such that
 \begin{align*}   gg'&=(g\otimes g')\phi  \\
                         g^{-1}  &=g\k,  \end{align*}
for each $g,g'\in G_{cl}$. The counit $\e\colon  \cal{A}\rightarrow
\Bbb{C}$  is
the neutral element of $G_{cl}$. We shall assume that the (complex)
Lie  algebra $\lie(G_{cl})$ is realized \cite{D} as  the
space  of  linear  functionals $X\colon \cal{A}\rightarrow \Bbb{C}$ satisfying
                      $$X(ab)=\e(a)X(b)+\e(b)X(a),$$
for each $a,b\in \cal{A}$.

Let $\Gamma$ be a first-order  differential   calculus   over
$G$. This means \cite{W3} that $\Gamma$ is a bimodule over  $\cal{A}$
endowed with a differential
$d\colon  \cal{A}\rightarrow  \Gamma$
such  that  elements  of   the   form
$a\,db$ linearly   generate    $\Gamma$.    Let
$$\Gamma^{\otimes} =\sideset{}{^\oplus}\sum_{k\geq 0} \Gamma^{\otimes
k}$$   be the tensor bundle
algebra \cite{W3} built  over
$\Gamma$. Let $$\Gamma^{\wedge} =\sideset{}{^\oplus}\sum_{k\geq 0}
\Gamma^{\wedge k}$$  be the universal
differential envelope (\cite{D}--Appendix B) of
$\Gamma$. The algebra $\Gamma^{\wedge}$  can  be  obtained
from  $\Gamma^{\otimes}$
by     factorising     through     the      ideal      $S^{\wedge}
\subseteq\Gamma^{\otimes}$   generated  by
the elements of  the   form
$$Q=\sum_ida_i  \otimes_{\cal{A}}  db_i,$$  where  $a_i ,b_i \in
\cal{A}$   satisfy  $\Sum_i a_i db_i =0$.   In    particular,
the differential $d\colon\Gamma^{\wedge}\rightarrow\Gamma^{\wedge}$
extends  $d\colon\cal{A}\rightarrow\Gamma,$ in a natural manner.

Let us assume that $\Gamma$  is  left-covariant  \cite{W3}  and  let
$\ell_{\Gamma} \colon \Gamma\rightarrow \cal{A}\otimes \Gamma$  be  the  left
action of $G$ on $\Gamma$. Let $\Gamma_{\inv}$    be  the  space  of
left-invariant elements of $\Gamma$ (playing the role  of  the  dual
space   of   the   Lie   algebra   of   $G$)   and   let   $\pi\colon
\cal{A}\rightarrow \Gamma_{\inv}$  be the canonical projection map, given by
  $$\pi(a)=\k(a^{(1)} )da^{(2)} .$$
This map is surjective  and $\cal{R}=\ker(\e)\cap  \ker(\pi)$  is  the
right  $\cal{A}$-ideal  which  canonically  \cite{W3}   corresponds   to
$\Gamma$.

 The   space   $\Gamma_{\inv}$      possesses   a    natural    right
$\cal{A}$-module structure, which will be denoted by $\circ $.
Explicitly
     $$\pi(a)\circ b=\pi\bigl[\bigl(a-\e(a)1\bigr)b\bigr]  $$
for each $a,b\in \cal{A}$.

 Let us  now  assume  that  $\Gamma$  is  bicovariant,  and  let
 $\wp_{\Gamma}\colon \Gamma\rightarrow \Gamma\otimes \cal{A} $
be  the  right
action of $G$ on $\Gamma$.

 The ``adjoint'' action $\ad\colon  \cal{A}\rightarrow  \cal{A}\otimes
\cal{A}$ of $G$ on $G$ is given by
                      $$\ad(a)=a^{(2)}\otimes \k(a^{(1)})a^{(3)}.$$

 The   space   $\Gamma_{\inv}$      is   right-invariant,   that   is
 $\wp_{\Gamma}(\Gamma_{\inv})\subseteq\Gamma_{\inv}   \otimes
\cal{A}$. The
corresponding  restriction   $\Ad\colon \Gamma_{\inv}
\rightarrow
\Gamma_{\inv}   \otimes \cal{A}$ is  interpretable  as  the  adjoint
action of $G$ on $\Gamma_{\inv}   $. Explicitly $\Ad$ is
characterized by
                  $$\Ad\pi=(\pi\otimes \id)\ad.$$

 The actions $\ell_{\Gamma}$  and  $\wp_{\Gamma}$ can  be  naturally
extended  to  the
grade  preserving  homomorphisms   $\wp_{\Gamma}^{\wedge,\otimes} \colon
\Gamma^{\wedge,\otimes} \rightarrow
\Gamma^{\wedge,\otimes} \otimes     \cal{A}$ and
$\ell_{\Gamma}^{\wedge,\otimes} \colon
\Gamma^{\wedge,\otimes} \rightarrow
\cal{A}\otimes \Gamma^{\wedge,\otimes}$ (their    restrictions    on
$\cal{A}$  coincide with $\phi$).

 The symbol $\grten$  will be used for the graded tensor product of
graded-differential  algebras.  The  comultiplication  $\phi$
admits the unique extension $\widehat{\phi}\colon \Gamma^{\wedge}\rightarrow
\Gamma^{\wedge} \grten  \Gamma^{\wedge}$
which is a homomorphism of graded-differential  algebras \cite{D}.  In
particular,
 $$ \widehat{\phi}(\xi)=\ell_{\Gamma}(\xi)+\wp_{\Gamma}(\xi)$$
for each $ \xi\in \Gamma$.
The antipode $\k$ admits the unique
extension $\widehat{\k}\colon\Gamma^\wedge\rightarrow
\Gamma^\wedge$, which is graded-antimultiplicative and satisfies
$\widehat{\k}d=d\widehat{\k}$.

     Let   us   denote    by    $\Gamma_{\inv}^{\otimes}$       and
$\Gamma_{\inv}^{\wedge}$    subalgebras of
left-invariant     elements     of     $\Gamma^{\otimes}$      and
$\Gamma^{\wedge}$  respectively. We have
$$\Gamma_{\inv}^{\otimes}=
\sideset{}{^\oplus}\sum_{k\geq 0} \Gamma_{\inv}^{\otimes k}\qquad
                          \Gamma_{\inv}^{\wedge}                 =
\sideset{}{^\oplus}\sum_{k\geq 0}  \Gamma_{\inv}^{\wedge  k}    ,
$$
where $\Gamma_{\inv}^{\otimes k}$    and $\Gamma_{\inv}^{\wedge k}$
consist of left-invariant
elements  from  $\Gamma^{\otimes k}$
and    $\Gamma^{\wedge    k}$      respectively.     The     space
$\Gamma_{\inv}^{\otimes k}$    is actually the
tensor product
of $k$-copies of $\Gamma_{\inv}$.

  The following natural
isomorphism holds
$$\Gamma_{\inv}^{\wedge}=\Gamma_{\inv}^{\otimes}/S_{\inv}^{\wedge},$$
where  $S_{\inv}^{\wedge}$ is  the  left-invariant   part   of
$S^{\wedge}$.  This  space  is  an
ideal  in  $\Gamma_{\inv}^{\otimes}$ generated   by   elements   of
the  form
$$q=\pi(a^{(1)}   )\otimes \pi(a^{(2)}   ),$$ where $a\in \cal{R}$.

 All  introduced  spaces   of   the   form   $\Gamma_{\inv}^*$      are
right-invariant.  We  shall  denote  by   $\Ad^{*}$ the  adjoint
actions of $G$ on the corresponding spaces.

     The formula
                     $$\vartheta\circ   a=\k(a^{(1)}     )\vartheta
a^{(2)} $$
defines an extension  of  the  right  $\cal{A}$-module   structure
$\circ$
from $\Gamma_{\inv}$    to $\Gamma_{\inv}^{\wedge,\otimes}$.  We
have
\begin{align*} 1\circ a&=\e(a)1 \\
                         (\vartheta\eta)\circ
a&=(\vartheta\circ a^{(1)})(\eta\circ a^{(2)}   )   \end{align*}
for each $\vartheta,\eta\in \Gamma_{\inv}^{\wedge,\otimes}$    and
$a\in \cal{A}$.

 The algebra $\Gamma_{\inv}^{\wedge}   \subseteq\Gamma^{\wedge}$  is
$d$-invariant.  The
differential    $d\colon    \Gamma_{\inv}^{\wedge}      \rightarrow
\Gamma_{\inv}^\wedge$    is
explicitly determined by
                   $$d\pi(a)=-\pi(a^{(1)})\pi(a^{(2)}).$$

 If  $\Gamma$  is  *-covariant  then  the  *-involution $ *\colon
\Gamma\rightarrow \Gamma$ is naturally extendible  from  $\Gamma$  to
$\Gamma^{\wedge ,\otimes}$    (such that for each $\vartheta,\eta\in
\Gamma^{\wedge ,\otimes}$    we  have
$(\vartheta\eta)^* =(-)^{\partial\vartheta\partial\eta}\eta^*  \vartheta^*
)$.           Algebras
$\Gamma_{\inv}^{\wedge}                     ,\Gamma_{\inv}^{\otimes}
\subseteq\Gamma^{\wedge ,\otimes}$    are *-invariant.
We have
              $$(\vartheta\circ a)^* =\vartheta^* \circ \k(a)^*$$
for each $a\in \cal{A}$ and $\vartheta\in \Gamma_{\inv}^{\wedge,\otimes}$.

     Explicitly, the *-involution on $\Gamma_{\inv}$    is determined by
                        $$   \pi(a)^* =-\pi[\k(a)^* ].$$

 The map $\widehat{\phi}$, as well as the left and the
right actions of $G$ on $\Gamma^{\wedge ,\otimes}$ are *-preserving,
in a natural manner.

 Let $M$ be a compact smooth manifold.
By definition \cite{D} a {\it quantum  principal
$G$-bundle}   over   $M$    is    a    triplet
$P=(\cal{B},i,F)$
where  $\cal{B}$  is  a  (unital) *-algebra,  consisting   of
appropriate ``functions'' on  $P$,  while  $F\colon  \cal{B}\rightarrow
\cal{B}\otimes \cal{A}$ and  $i\colon  S(M)\rightarrow  \cal{B}$  are
(unital) *-homomorphisms,  intrepretable  as  the  dualized  right
action of $G$ on $P,$ and the dualized  projection   of   $P$   on
$M$. Further, the
bundle $P$ is {\it locally trivial} in the sense that for each  $x\in
M$ there exists an open set $U\subseteq M$  such  that  $x\in  U,$  and  a
*-homomorphism  $\pi_U  \colon   \cal{B}\rightarrow    S(U)\otimes
\cal{A}$
such that
\begin{gather*}
\pi_U i(f)=(f{\restr}_U )\otimes 1\\
\pi_U (\cal{B})\supseteq \SC (U)\otimes \cal{A}\\
(\id\otimes \phi)\pi_U =(\pi_U \otimes \id)F,
\end{gather*}
and such that
$$\pi_U\bigl(i(f)b\bigr)=0\quad\Longrightarrow\quad i(f)b=0$$
for each $f\in  \SC  (U)$. Here $S$ and $\SC$ denote the corresponding
*-algebras of complex smooth functions (with compact supports, respectively).

The  homomorphism  $\pi_U$   is   interpretable   as   the   dualized
trivialization of $P$ over $U$. Every pair $(U,\pi_U ),$ consisting of  an
open  set  $U\subseteq M$  and   of   a   *-homomorphism   $\pi_U \colon
\cal{B}\rightarrow   S(U)\otimes   \cal{A},$    satisfying    above
conditions is called a local trivialization of $P$.

A trivialization system  for  $P$  is  a  family
$\tau=\bigl\{(U,\pi_U  )\mid U\in\cal{U}\bigr\}$  of
local trivializations of $P$, where $\cal{U}$ is a finite  open  cover
of $M$.

For each $k\in \Bbb{N}$ we shall denote by $N^k (\cal{U})$
the  set  of
$k$-tuples $(U_1 ,\dots ,U_k )\in \cal{U}^k$  such that $U_1 \cap
\dots \cap U_k\neq \emptyset$.

     The  main  structural  result  concerning  quantum  principal
bundles is that there  exists  a  natural  correspondence  between
quantum  principal  $G$-bundles  $P$   and   classical   principal
$G_{cl}$-bundles
$P_{cl}$   over $M.$ This corresponence can be described as follows.

 From a  given  trivialization  system  $\tau$  it  is  possible  to
construct the corresponding  $G$-cocycle  which  is  a  system  of
*-automorphisms    $\psi_{UV}$ of $S(U\cap V)\otimes
\cal{A}$, where $(U,V)\in   N^2
(\cal{U})$,
realizing transformations  between  $(V,\pi_V )$  and
$(U,\pi_U )$. Such systems of maps completely determine the bundle
$P$.

     Explicitly, let us consider the *-algebra
               $$\Sigma  (\cal{U})=  \sideset{}{^\oplus}\sum_{U\in
\cal{U}} \Bigl[S(U)\otimes \cal{A}\Bigr].$$
The algebra $\cal{B}$ is realizable as a subalgebra  of  $\Sigma
(\cal{U}),$ consisting of elements $b\in \Sigma (\cal{U})$ satisfying
        $$(_U {\restr}_{U\cap V} \otimes \id)p_U (b)=\psi_{UV}
(_V {\restr}_{U\cap V}  \otimes \id)p_V (b) $$
for   each   $(U,V)\in N^2 (\cal{U}).$   Here,    $p_U \colon    \Sigma
(\cal{U})\rightarrow   S(U)\otimes    \cal{A}$    are    coordinate
projections. In terms of this realization we have
                           $$ \pi_U =p_U{\restr}\cal{B},$$
for each $U\in \cal{U}.$

 However, it turns out that $G$-cocycles are in a natural  bijection
with standard $G_{cl}$-cocycles (over $\cal{U}$), which  are  systems  of
smooth maps $g_{UV}  \colon U\cap V\rightarrow G_{cl}$  satisfying
$$g_{UV}g_{VW}(x)=g_{UW}(x),$$
for each $(U,V,W)\in N^3(\cal{U})$ and $x\in U\cap V\cap W$
(in particular $g_{UV}^{-1}=g_{VU}$). The correspondence is established
via the following formula
$$\psi_{UV}(\varphi\otimes a)=\varphi g_{VU}(a^{(1)})\otimes a^{(2)}.$$
Here, maps $g_{UV}$ are understood as *-homomorphisms $g_{UV}\colon
\cal{A}\rightarrow S(U)$, in a natural manner.
On the other hand, $G_{cl}$-cocycles determine, in the standard manner,
principal $G_{cl}$-bundles $P$ over $M.$

 The bundle $P_{cl}$   is interpretable as the ``classical part'' of
$P.$ The
elements of $P_{cl}$   are in a natural bijection  with  *-characters  of
$\cal{B}.$  The  correspondence  $P\leftrightarrow  P_{cl}$    has  a   simple
geometrical    explanation.    The     ``transition      functions''
$\psi_{UV}$
are, at  the  geometrical  level,  vertical  ``diffeomorphisms''  of
$(U\cap V)\times G.$ Therefore they preserve the geometrical structure  of
$(U\cap V)\times G.$ In particular, they must preserve the classical  part
$(U\cap V)\times G_{cl}$   consisting  of  points  of  $(U\cap
V)\times G.$  Moreover,
transition diffeomorphisms  are  completely  determined  by  their
restrictions on $(U\cap V)\times G_{cl}  ,$ because of  the  right
covariance.
The corresponding ``restrictions'' are precisely transition  functions
for the classical bundle $P_{cl}.$

 So far about the structure of quantum principal bundles. For each
(nonempty) open  set  $U\subseteq M$  let  $\Omega (U)$  be  the
graded-differential *-algebra  of differential forms on $U.$ In  developing
a differential calculus  over  quantum  principal  bundles  it  is
natural to assume that the calculus is fully compatible  with
the geometrical structure on  the  bundle,  such  that  all  local
trivializations of the bundle locally trivialize the calculus  too
(a precise formulation of this condition is given in
\cite{D}--Section~3). It turns out that this  condition  completely  fixes  the
calculus on the bundle (if the calculus on the  structure  quantum
group  is  fixed).  However,   the   condition   implies   certain
restrictions on a possible differential calculus  $\Gamma$  over $G.$

Namely,  all  retrivialization  maps  $\psi_{UV}$ must   be
extendible     to     differential      algebra automorphisms
$\psi^{\wedge}_{UV}\colon\Omega(U\cap V)\grten \Gamma^{\wedge}\rightarrow
\Omega(U\cap  V)\grten
\Gamma^{\wedge} .$ Differential calculi $\Gamma$ satisfying this condition
are called admissible. If $\Gamma$ is left-covariant then it  is
admissible  iff  $$(X\otimes  \id)\ad(\cal{R})=\{0\}$$   for   each   $X\in
\lie(G_{cl}  ).$  This  fact  implies  that  there  exists  the   minimal
admissible left-covariant calculus $\Gamma.$  This  calculus  is
based         on         the         right         $\cal{A}$-ideal
$\widehat{\cal{R}}\subseteq \ker(\e)$ consisting
of all elements $a\in \ker(\e)$ satisfying $$(X\otimes  \id)\ad(a)=0,$$  for
each      $X\in       lie(G_{cl}  ). $

Moreover, we have
$\ad(\widehat{\cal{R}})\subseteq\widehat{\cal{R}}\otimes \cal{A}$ and
$\k(\widehat{\cal{R}})^* =\widehat{\cal{R}}$,
which implies \cite{W3} that $\Gamma$ is bicovariant and *-covariant
respectively.

 In the following, $\Gamma$  will  be  this  minimal  admissible
(bicovariant *-) calculus.

 Let $\Omega(P)$ be the graded-differential  *-algebra  representing
differential calculus on $P$ (constructed by combining  differential
forms on $M$ with the universal  envelope  $\Gamma^{\wedge}$   of
$\Gamma).$
Explicitly, let us consider the direct sum
         $$\Sigma^{\wedge}(\cal{U})=
\sideset{}{^\oplus}\sum_{U\in\cal{U}} \Bigl[\Omega(U)\grten
\Gamma^{\wedge} \Bigr].$$
Then $\Omega(P)$ can be viewed as a graded-differential  subalgebra
consisting of elements $w\in \Sigma^{\wedge}  (\cal{U})$ satisfying
              $$(_U{\restr}_{U\cap V}\otimes\id)p_U
(w)=\psi^{\wedge}_{UV}(_V{\restr}_{U\cap V}\otimes\id)p_V (w) $$
for    each      $(U,V)\in      N^2  (\cal{U})$. Here $p_U  \colon
\Sigma^{\wedge}
 (\cal{U})\rightarrow  \Omega(U)\grten\Gamma^{\wedge}$   are
corresponding coordinate projections.

 As  a   differential   algebra,   $\Omega(P)$   is   generated   by
$\cal{B}=\Omega^0 (P).$ For every local  trivialisation  $(U,\pi_U
)$  of  $P$
there  exists  the  unique   differential   algebra   homomorphism
$\pi^{\wedge}_U            \colon             \Omega(P)\rightarrow
\Omega(U)\grten \Gamma^{\wedge}$  extending
$\pi_U$    (in   fact   $\pi^{\wedge}_U =p_U {\restr}\Omega(P)). $  The
map    $i\colon
S(M)\rightarrow    \cal{B}$    admits    a    natural    extension
$i^{\wedge} \colon
\Omega(M)\rightarrow \Omega(P),$  which  is  interpretable  as  the
``pull back'' of differential forms on $M$ to $P$. We have
$$\pi^{\wedge}_Ui^{\wedge}(w)=(w{\restr}_U)\otimes 1.$$

     The right  action  $F\colon \cal{B}\rightarrow \cal{B}\otimes
\cal{A}$
is  (uniquely) extendible  to  a
differential algebra homomorphism $\widehat{F}\colon \Omega(P)\rightarrow
\Omega(P)\grten \Gamma^{\wedge}$,
imitating the corresponding pull back map. The formula
                         $$     F^{\wedge}     =(\id\otimes \Pi
)\widehat{F}$$
determines   a   *-homomorphism   $F^{\wedge} \colon    \Omega(P)\rightarrow
\Omega(P)\otimes \cal{A}$ interpretable  as  the  (dualized)  right
action of $G$ on $\Omega(P)$. Here
$\Pi\colon\Gamma^\wedge\rightarrow\cal{A}$ is the projection map.

 Let $\ver(P)$ be the  graded-differential  *-algebra  obtained  by
factorizing $\Omega(P)$ through the (differential *-ideal) generated
by elements of the form $d[i(f)]$. The  elements  of
$\ver(P)$ play the role of ``verticalized'' differential   forms  on
$P$ (in classical geometry, entities obtained by restricting the domain
of differential forms to the Lie algebra of vertical vector fields
on the bundle).
At the level of graded  vector  spaces,  there  exists  a  natural
isomorphism
         $$\ver(P)\cong \cal{B}\otimes \Gamma^{\wedge}_{\inv}   .$$

Let $\pi_v \colon \Omega(P)\rightarrow \ver(P)$ be  the  corresponding
projection  map.  In  terms  of  the  above  identifications,  the
differential *-algebra structure on $\ver(P)$ is specified by
\begin{gather*}
 (q\otimes  \eta)(b\otimes  \vartheta)=\sum_k qb_k \otimes   (\eta\circ
a_k )\vartheta  \\
 (b\otimes \vartheta)^* =\sum_k b^*_k  \otimes  (\vartheta^*  \circ
a^*_k)\\
d_v  (b\otimes   \vartheta)=\sum_k  b_k  \otimes
\pi(a_k )\vartheta+b\otimes
d\vartheta
\end{gather*}
where $F(b)=\Sum_k b_k\otimes a_k$.

 Another important algebra naturally associated to $\Omega(P)$ is  a
graded   *-subalgebra    $\hor(P)\subseteq\Omega(P) $   representing
horizontal forms. By definition, $\hor(P)$ consists  of  forms  $w\in
\Omega(P)$ with the property $$\pi^{\wedge}_U (w)\in
\Omega(U)\otimes  \cal{A},$$
for  each  local  trivialization  $(U,\pi_U  )$. Equivalently,
$$\hor(P)=(\widehat{F})^{-1}\Bigl\{\Omega(P)\otimes\cal{A}\Bigr\}.$$
The    algebra
$\hor(P)$  is
invariant  under  the  right  action  of   $G,$   in   other   words
$$F^{\wedge} \bigl(\hor(P)\bigr)\subseteq\hor(P)\otimes \cal{A}.$$

 Let $\psi(P)$ be  the  space  of  all  linear  maps  $\varphi\colon
\Gamma_{\inv}   \rightarrow \Omega(P)$ satisfying
              $$   (\varphi\otimes    \id)\Ad=F^{\wedge}\varphi. $$
This space is  naturally  graded  (the  grading  is  induced  from
$\Omega(P)$). The elements of $\psi (P)$ are quantum  counterparts  of
pseudotensorial forms on  the  bundle  with  coefficients  in  the
structure   group   Lie   algebra   (relative   to   the   adjoint
representation). The space $\psi (P) $ is  closed  with  respect  to
compositions with $d\colon \Omega(P)\rightarrow \Omega(P).$

 Let $\tau (P)\subseteq\psi (P)$ be the subspace consisting  of
$\hor(P)$-valued maps. This space is imaginable as consisting  of
the corresponding tensorial forms.

There exists a natural *-involution on $\psi (P)$. It is given by
           $$\varphi^* (\vartheta)=\varphi(\vartheta^* )^* .$$
The space $\tau (P)$ is *-invariant.

 Tensorial forms possess the following local representation:
     $$\pi^{\wedge}_U       \varphi(\vartheta)=(f^U        \otimes
\id)\Ad(\vartheta),$$
where $f^U \colon\Gamma_{\inv} \rightarrow \Omega(U)$ is a linear map.

 For the purposes of this paper the most important  topic  of  the
theory  of  quantum  principal  bundles  is   the   formalism   of
connections.  By  definition,  a  connection   on   $P$   is   every
pseudotensorial hermitian $1$-form $\omega$ satisfying
     $$\pi_v \omega(\vartheta)=1\otimes \vartheta  $$
for each $\vartheta\in \Gamma_{\inv}$. The above  formula  is  the
quantum counterpart for the classical condition  that  connections
map fundamental vector fields into their generators. Connections form a
real affine space $\con(P)$.

     In  local  terms,  connections possess  the  following
representation
      $$\pi^{\wedge}_U \omega(\vartheta)=(A^U \otimes
\id)\Ad(\vartheta)+1_U \otimes
\vartheta, $$
where $A^U \colon \Gamma_{\inv}   \rightarrow \Omega(U)$ is a  $1$-form
valued hermitian linear map (playing the role of the corresponding
gauge potential).

     The curvature operator can be described as follows.

 Let  us  fix  a   map   $\delta\colon   \Gamma_{\inv}   \rightarrow
\Gamma_{\inv}   \otimes  \Gamma_{\inv}$     which   intertwines   the
corresponding adjoint actions and such that if
 $$\delta(\vartheta)=\sum_k \vartheta^1_k \otimes \vartheta^2_k $$
then
$$\delta(\vartheta^* )=-\sum_k (\vartheta^2_k )^* \otimes
(\vartheta^1_k )^* \qquad
d\vartheta=\sum_k\vartheta^1_k\vartheta^2_k.$$
Every such a map will be called an {\it embedded differential}. Further,
for   each   pair   of    linear    maps    $\varphi,\psi$ on $
\Gamma_{\inv}$ with values in an  arbitrary
algebra $\Omega$ let  $\langle \varphi,\psi\rangle\colon
\Gamma_{\inv}   \rightarrow
\Omega$ be a map given by
$$\langle \varphi,\psi\rangle(\vartheta)=\sum_k\varphi(\vartheta_k^1)
\psi(\vartheta_k^2).$$
By construction, if $\varphi,\psi \in \psi (P)$ then
$\langle \varphi,\psi\rangle\in  \psi
(P),$ too.

Finally, the curvature $R_\omega$   of a connection $\omega $ can
                       be       defined        as
$$  R_\omega =d\omega-\langle \omega,\omega\rangle  . $$
The above formula corresponds to the  structure  equation  in  the
classical theory. It turns out that $R_\omega$    is  a  tensorial
$2$-form. Locally, in terms of the corresponding gauge potentials  we
have
        $$\pi^{\wedge}_U R_\omega  (\vartheta)=(F^U \otimes
\id)\Ad(\vartheta),$$
 where
                    $$  F^U =dA^U -\langle A^U ,A^U \rangle  . $$

 For each open set $U\subseteq M$ the symbol $\otimes_U$   will  be  used
for the tensor  product  over  $S(U)$.
Similarly, the symbol $\grten_U$ will denote  the
graded tensor product of graded-differential *-algebras containing
$\Omega(U)$ as their subalgebra.

\section{ Quantum Gauge Bundles }
     This section  is  devoted  to  generalizations  of  the  most
important aspects of the concept of gauge transformations, in  the
framework of the formalism of quantum principal bundles. The  main
geometrical object that will be constructed is the  quantum  gauge
bundle, a noncommutative-geometric counterpart of the gauge bundle
of the classical theory.
\subsection{Classical Consideration}
     In order to present  motivations  for  contructions  of  this
section let us assume for a moment that $G$ is an  ordinary  compact
Lie group, and let $P$ be a (classical) principal bundle over $M.$

 By  definition,  gauge  transformations   of   $P$   are   vertical
automorphisms   of   this   bundle.   In   other   words,    gauge
transformations are diffeomorphisms  $\psi   \colon   P\rightarrow
P$
satisfying
\begin{align*}   &\pi_M\psi =\pi_M \\
&\psi (pg)=\psi (p)g, \end{align*}
for each $p\in P$ and $g\in G,$ where $(p,g)\mapsto pg$ is the right
action of $G$ on $P$ and $\pi_M\colon P\rightarrow M$  is  the  projection
map. Equivalently,  gauge  transformations  are  interpretable  as
(smooth) sections of the gauge bundle $\ADP,$ which is  the  bundle
associated to $P,$ with respect to the  adjoint  action  of  $G$  onto
itself.

     The equivalence between two definitions  is  established  via
the following formula
                         $$ \psi (p)=pf(p), $$
where $f\colon P\rightarrow G$ is a smooth equivariant function in
the sense that
                       $$   f(pg)=g^{-1}  f(p)g, $$
for each $p\in P$ and $g\in  G.$  Such  functions  are  in  a
natural  correspondence  with  sections  of  the  corresponding  associated
bundle $\ADP.$

 For each $x\in M$ the fiber $G_x =\pi_M^{\sharp-1}  (x)$ over $x$
(where  $\pi_M^\sharp\colon
\ADP\rightarrow M$ is the projection map) possesses a natural  Lie
group  structure.  The   group   $G_x$    is   isomorphic   (generally
non-invariantly) to $G.$ For a given $p\in \pi_M^{-1}  (x)=P_x$
there exists  a canonical diffeomorphism $G\leftrightarrow P_x$
defined by $g\leftrightarrow pg$, and a group isomorphism
$G\leftrightarrow    G_x  $
given by $g\leftrightarrow \bigl[(p,g)\bigr]$.
Here, $\ADP$ is understood  as  the  orbit
space of the right action $((p,g'),g)\mapsto (pg,g^{-1}  g'g)$
of  $G$
on $P\times G$ and $\bigl[\,\bigr]$ denotes the corresponding orbit.

     There exists a natural left action of $G_x$  on $P_x .$ In  terms  of
the above identifications this action becomes the
multiplication  on the left. Collecting all these fiber actions together,
we obtain a smooth map
 \begin{equation}\label{33}\beta^*_M \colon \ADP\times_M
P\rightarrow P. \end{equation}
With the help of $\beta^*_M$  the equivalence between gauge  transformations
$\psi$   and  sections  $\varphi\colon  M\rightarrow  \ADP$  can   be
described as follows
    \begin{equation}\label{34}\psi   =\beta^*_M   (\varphi\times_M
\id). \end{equation}

 Moreover, the correspondence  $\psi  \leftrightarrow  \varphi$  is  an
isomorphism between the group $\cal{G}$ of gauge transformations of $P,$ and
the group $\Gamma(\ADP)$ of smooth sections of $\ADP.$

     The  group  structure  in  fibers  of  $\ADP$  determine   the
following maps of bundles
\begin{equation}\label{35}\begin{gathered}
\text{\it the fibewise multiplication}\quad
 \phi^*_M \colon \ADP\times_M \ADP\rightarrow \ADP\\
 \text{\it  the   unit section}\qquad\qquad\quad \e^*_M \colon M\rightarrow
\ADP\\
\text{\it the fiberwise inverse}\qquad\quad \k^*_M\colon \ADP
\rightarrow \ADP\end{gathered} \end{equation}

     At the dual level of function algebras \eqref{33}  and  \eqref{35}  are
represented by the corresponding $S(M)$-linear *-homomorphisms
\begin{equation}\label{36}\begin{aligned}
\phi_M &\colon S(\ADP)\rightarrow S(\ADP)\otimes_M  S(\ADP)\\
\e_M &\colon S(\ADP)\rightarrow S(M)  \\
\k_M &\colon S(\ADP)\rightarrow S(\ADP)  \\
\beta_M &\colon S(P)\rightarrow S(\ADP)\otimes_M  S(\ADP).
\end{aligned} \end{equation}

 The structure of the gauge group is completely  encoded  in  maps
$\bigl\{\phi_M,\k_M,\e_M\bigr\}$.

 At the dual level, gauge transformations $\psi$  can be  viewed  as
$S(M)$-linear *-automorphisms $\psi\colon S(P)\rightarrow S(P)$
intertwining  the  (dualized)  right   action   of   $G.$   Further,
interpreted as sections of $\ADP,$ gauge transformations become, at
the  dual  level,   $S(M)$-linear   *-homomorphisms   $\varphi \colon
S(\ADP)\rightarrow S(M).$ In this picture,  the  action  of  $\cal{G}$  on
$S(\ADP)$ is given by
           $$(\varphi,f)\mapsto(\varphi\otimes \id)\beta_M (f).$$

 The maps \eqref{36} are not suitable  for  considering  situations  in
which gauge transformations act on differential forms. This can be
easily ``improved'' by  extending  these  maps  to  $\Omega(M)$-linear
homomorphisms
\begin{equation}\label{37}\begin{aligned}
   \widehat{\phi}_M &\colon \Omega(\ADP)\rightarrow
    \Omega(\ADP)\grten_M
                   \Omega(\ADP)  \\
   \widehat{\e}_M &\colon \Omega(\ADP)\rightarrow \Omega(M) \\
  \widehat{\k}_M &\colon \Omega(\ADP)\rightarrow \Omega(\ADP) \\
 \widehat{\beta}_M &\colon \Omega(P)\rightarrow \Omega(\ADP)
 \grten_M  \Omega(\ADP) \end{aligned} \end{equation}
of graded-differential *-algebras. It is worth noticing  that  the
above maps are unique, as graded-differential extensions. Actually these maps
can be viewed as ``pull backs'' of \eqref{33}
 and \eqref{35}.
\subsection{Quantum Consideration}

The presented picture admits a direct noncommutative-geometric generalization.
As first, we shall construct,  starting
from a quantum principal bundle $P$, the corresponding quantum gauge
bundle  $\ADP$.  Then  the  counterparts  of  maps \eqref{36}  will  be
introduced and analyzed. In analogy with  the  classical  case  we
shall define gauge transformations as  vertical  automorphisms  of
the  bundle  $P$. It  turns  out  that  such   gauge
transformations of $P$ are in  a  natural  bijection  with  ordinary
gauge transformations of the classical part $P_{cl}$ of  $P$.  We  shall
also   study   various   equivalent   interpretations   of   gauge
transformations. Finally, a canonical differential calculus on the
bundle $\ADP$ will be constructed and analyzed.

 Let $G$ be a compact matrix quantum group, and let  $P=(\cal{B},i,F)$
be  a  quantum  principal  $G$-bundle  over  $M$.   Let   us   fix   a
trivialization system $\tau$ for $P$. For  each  $(U,V)\in  N^2 (\cal{U})$
let  us  define  a
linear   map    $\xi_{UV}\colon S(U\cap V)\otimes\cal{A}\rightarrow
S(U\cap V)\otimes \cal{A}$ by the following formula
\begin{equation}\label{38} \xi_{UV}(\varphi\otimes a)
=\varphi g_{UV}\bigl[\k(a^{(1)})a^{(3)}\bigr]\otimes a^{(2)}.
\end{equation}

\begin{lem}\label{lem:31}
(i) The maps $\xi_{UV}$ are  $S(U\cap V)$-linear  *-automorphisms
and
\begin{equation}\xi^{-1}_{UV}=\xi_{VU}.\label{xinv}\end{equation}

 (ii)  We have
\begin{equation}\label{39}\xi_{UV}\xi_{VW}(\varphi)=\xi_{UW}(\varphi).
\end{equation}
for  each  $(U,V,W)\in  N^3 (\cal{U})$  and $\varphi\in\SC(U{\cap}
V{\cap} W)\otimes\cal{A}$.

    (iii) The diagrams
\begin{gather*}
\begin{CD}
S(U\cap V)\otimes\cal{A}@>{\mbox{$\id\otimes\phi$}}>>
\bigl[S(U\cap V)\otimes\cal{A}
\bigr]\otimes_{U\cap V}\bigl[S(U\cap V)\otimes\cal{A}\bigr]\\
@V{\mbox{$\xi_{UV}$}}VV @VV{\mbox{$\xi_{UV}\otimes
\xi_{UV}$}}V\\
S(U\cap V)\otimes\cal{A} @>>{\mbox{$\id\otimes\phi$}}>
\bigl[S(U\cap V)\otimes\cal{A}
\bigr]\otimes_{U\cap V}\bigl[S(U\cap V)\otimes\cal{A}\bigr]
\end{CD}\\
\begin{CD}
S(U\cap V)\otimes\cal{A}@>{\mbox{$\id\otimes \e$}}>> S(U\cap V)\\
@V{\mbox{$\xi_{UV}$}}VV @VV{\mbox{$\id$}}V\\
S(U\cap V)\otimes\cal{A}@>>{\mbox{$\id\otimes \e$}}> S(U\cap V)
\end{CD}\qquad
\begin{CD}
S(U\cap V)\otimes\cal{A}@>{\mbox{$\id\otimes \k$}}>> S(U\cap V)\otimes\cal{A}\\
@V{\mbox{$\xi_{UV}$}}VV @VV{\mbox{$\xi_{UV}$}}V\\
S(U\cap V)\otimes\cal{A}@>>{\mbox{$\id\otimes \k$}}> S(U\cap V)\otimes\cal{A}
\end{CD}\\
\begin{CD}
S(U\cap V)\otimes\cal{A}@>{\mbox{$\id\otimes\phi$}}>>
\bigl[S(U\cap V)\otimes\cal{A}
\bigr]\otimes_{U\cap V}\bigl[S(U\cap V)\otimes\cal{A}\bigr]\\
@V{\mbox{$\psi_{UV}$}}VV @VV{\mbox{$\xi_{UV}\otimes
\psi_{UV}$}}V\\
S(U\cap V)\otimes\cal{A} @>>{\mbox{$\id\otimes\phi$}}>
\bigl[S(U\cap V)\otimes\cal{A}
\bigr]\otimes_{U\cap V}\bigl[S(U\cap V)\otimes\cal{A}\bigr]
\end{CD}
\end{gather*}
are commutative. \end{lem}

\begin{pf}
We have
\begin{equation*} \begin{split}
\xi_{UV}    \xi_{VW}    (\varphi\otimes
a)&=\xi_{UV} \left(\varphi g_{VW}\bigl[\k(a^{(1)}   )a^{(3)}
\bigr]\otimes a^{(2)} \right)\\
&=\varphi g_{VW}\bigl[\k(a^{(1)}   )a^{(5)}\bigr]g_{UV}
\bigl[\k(a^{(2)})a^{(4)}\bigr]\otimes a^{(3)}  \\
&=\varphi g_{UW}\bigl[\k(a^{(1)})a^{(3)}\bigr]\otimes a^{(2)}
=\xi_{UW}    (\varphi\otimes a),\end{split}\end{equation*}
for  each  $(U,V,W)\in  N^3 (\cal{U})$,
$ \varphi\in \SC(U{\cap} V{\cap} W)$ and $a\in \cal{A}$.
In particular, for $W=V$ this  implies
that the maps $\xi_{UV}$
are bijective and that \eqref{xinv} holds.

     The maps $\xi_{UV}$     are *-homomorphisms because of
\begin{multline*}
\xi_{UV}  (\varphi^* \otimes a^* )=\varphi^* g_{VU}
(a^{(1)*})g_{UV}(a^{(3)*}    )\otimes a^{(2)*}\\
=\left[\varphi g_{VU}(a^{(1)}   )g_{UV}(a^{(3)}
)\right]^* \!\otimes a^{(2)*}=\xi_{UV}(\varphi\otimes a)^*,
\end{multline*} and
\begin{equation*} \begin{split}
\xi_{UV}    (\varphi\psi
\otimes    ab)&=\varphi\psi
g_{VU}(a^{(1)}   b^{(1)}   )g_{UV}(a^{(3)}   b^{(3)}
   )\otimes a^{(2)}   b^{(2)}   \\
&=\left[\varphi g_{VU}(a^{(1)}   )g_{UV}(a^{(3)}   )\otimes
a^{(2)}   \right]\left[\psi
g_{VU}(b^{(1)}   )g_{UV}  (b^{(3)}   )\otimes b^{(2)}
 \right]\\
&=\xi_{UV}(\varphi\otimes a)\xi_{UV}(\psi \otimes b).
\end{split}\end{equation*}

     Finally, let us check commutativity of the above diagrams. We
compute
\begin{equation*}\begin{split}   (\xi_{UV}    \otimes \xi_{UV}    )
      (\id\otimes \phi)(\varphi\otimes a)
&=\varphi g_{UV}  \left(\k(a^{(1)}   )a^{(3)}\k(a^{(4)}
 )a^{(6)} \right)\otimes a^{(2)}   \otimes a^{(5)}  \\
&=\varphi g_{UV}  \left(\k(a^{(1)}   )a^{(4)} \right)\otimes
a^{(2)} \otimes a^{(3)}\\
&=(\id\otimes
\phi)\xi_{UV}    (\varphi\otimes a),\\
(\xi_{UV}    \otimes \psi_{UV}     )
 (\id\otimes \phi)(\varphi\otimes
a)&=\varphi g_{UV}  \bigl(\k(a^{(1)}   )a^{(3)}   \bigr)
 g_{VU}  (a^{(4)}   )
 \otimes a^{(2)}   \otimes a^{(5)}   \\
&=\varphi g_{VU}  (a^{(1)}   )\otimes a^{(2)}
\otimes a^{(3)}=(\id\otimes \phi)\psi_{UV}(\varphi\otimes a),\\
(\id\otimes \e)\xi_{UV}(\varphi\otimes
a)&=\varphi g_{UV}\left(\k(a^{(1)})a^{(2)}\right)=\e(a)\varphi,\\
\xi_{UV}    \bigl(\varphi\otimes \k(a)\bigr)&=\varphi g_{UV}
\bigl(\k^2 (a^{(3)}   )\k(a^{(1)}   )\bigr)\otimes \k(a^{(2)}   )  \\
&=\varphi g_{UV}  \bigl(\k(a^{(1)}   )a^{(3)}   \bigr)\otimes
               \k(a^{(2)}   )\\&=(\id\otimes
\k)\xi_{UV}    (\varphi\otimes a). \qed
\end{split}
\end{equation*}
\renewcommand{\qed}{}
\end{pf}

 The (algebra of functions on the) quantum gauge bundle $\ADP$  can
be now constructed as follows. Let $\cal{D}$ be the set of  elements
$q\in \Sigma (\cal{U})$ such that
 \begin{equation}\label{314}
  ({}_U {\restr}_{U\cap V}   \otimes \id)p_U (q)
   =\xi_{UV} ({}_V {\restr}_{U\cap V}   \otimes \id)p_V (q)
   \end{equation}
for each $(U,V)\in N^2 (\cal{U})$.

 Clearly, $\cal{D}$ is  a  *-subalgebra  of  $\Sigma  (\cal{U})$.  The
quantum space $\ADP$ corresponding to $\cal{D}$ plays the  role  of
the bundle associated to the principal bundle $P,$ with  respect  to
the adoint action  of  $G$  onto  itself  (represented  by  $\ad\colon
\cal{A}\rightarrow \cal{A}\otimes \cal{A}$). The fact that $\ADP$ is
a bundle over  $M$  is  established  throught  the  existence  of  a
*-monomorphism $j_M \colon S(M)\rightarrow \cal{D},$ playing the  role
of the dualized fibering of $\ADP$ over $M$. This map is defined by
equalities
\begin{equation}\label{315}
               p_U j_M (f)=(f{\restr} U)\otimes 1.
\end{equation}

\begin{defn}\label{defn:31}
The pair $\ADP=(\cal{D},j_M )$ is called {\it the quantum
gauge bundle} associated to $P$.
\end{defn}
 We  are  going  to  introduce  quantum   counterparts   of   maps
$\phi_M$, $\k_M$, $\e_M$   and  $\beta_M $.
For  each  $U\in  \cal{U}$, let $\pi^\sharp_U\colon
\cal{D}\rightarrow S(U)\otimes \cal{A}$ be the restriction of $p_U$  on
$\cal{D}$.

 \begin{pro}\label{pro:32}
 (i)  There  exist  the  unique   linear   maps
$\phi_M \colon\cal{D}\rightarrow   \cal{D}\otimes_M     \cal{D},$
$\e_M \colon
\cal{D}\rightarrow S(M)$, $\k_M \colon \cal{D}\rightarrow  \cal{D}$  and
$\beta_M \colon \cal{B}\rightarrow \cal{D}\otimes_M  \cal{B}$ such that
\begin{gather}
(\pi^\sharp_U\otimes \pi^\sharp_U)\phi_M
=(\id\otimes \phi)\pi^\sharp_U
\label{316}\\
(\pi^\sharp_U\otimes\pi_U)\beta_M
=(\id\otimes\phi)\pi^\sharp_U
\label{317}\\
\pi^\sharp_U \k_M =(\id\otimes \k)\pi^\sharp_U \label{318}\\
{\restr}_U \e_M =(\id\otimes \e)\pi^\sharp_U, \label{319}
\end{gather}
for each $U\in  \cal{U}$.  Here,
$S(U)\otimes   \cal{A}\otimes   \cal{A}$   and
$\bigl(S(U)\otimes\cal{A}\bigr)\otimes_U\bigl(S(U)\otimes\cal{A}\bigr)$
are identified, in a natural manner.

 (ii) All maps are $S(M)$-linear.  The  maps  $\phi_M$ ,$\e_M$   and
$\beta_M$   are
*-homomorphisms while $\k_M$  is antimultiplicative and
\begin{equation}\label{320}
\k_M \bigl[\k_M (f^* )^* \bigr]=f
\end{equation}
for each $f\in \cal{D}$.
\end{pro}
 \begin{pf} The  above  equalities  uniquely  fix  the  values  of  maps
$\phi_M$, $\e_M$, $\k_M$  and $\beta_M$
because  the  maps $\pi_U$ and
$\pi^\sharp_U$ distinguish points of $\cal{B}$ and $\cal{D}$.

 Let us  consider  the  algebra
$$\Sigma^*(\cal{U})= \sideset{}{^\oplus}\sum_{U\in \cal{U}} S(U)\otimes
\cal{A}\otimes  \cal{A}.$$
Algebras  $\cal{D}\otimes_M   \cal{D}$   and
$\cal{D}\otimes_M \cal{B}$   are understandable    as
subalgebras of $\Sigma^*(\cal{U})$. Let us consider maps
$\phi_M \colon
\Sigma (\cal{U})\rightarrow  \Sigma^*(\cal{U})$,
$\k_M\colon\Sigma(\cal{U})\rightarrow  \Sigma   (\cal{U})$   and
$\e_M \colon\Sigma(\cal{U})\rightarrow S(\cal{U})$ defined by
\begin{gather*}
p^*_U\phi_M =(\id\otimes \phi)p_U\\
p_U\k_M =(\id\otimes \k)p_U\\
{\restr}_U\e_M =(\id\otimes \e)p_U ,
\end{gather*}
where  $p^*_U\colon  \Sigma^*(\cal{U})\rightarrow
S(U)\otimes  \cal{A}\otimes
\cal{A}$ are coordinate projections and $S(\cal{U})$ is the direct sum of
algebras $S(U)$.

It is  easy  to  see  that
 $\phi_M (\cal{B})\subseteq\cal{D}\otimes_M
\cal{B}$, $\phi_M (\cal{D})\subseteq\cal{D}\otimes_M\cal{D}$,
$\k_M (\cal{D})\subseteq\cal{D}$  and
$\e_M (\cal{D})\subseteq S(M)$.   Let
us denote by $\bigl\{\phi_M, \beta_M, \k_M, \e_M\bigr\}$
the corresponding restrictions.
By construction,\eqref{316}--\eqref{319}
hold,
maps $\beta_M$ ,$\phi_M$  and $\e_M$  are *-homomorphisms,
$\k_M$  is antimultiplicative and \eqref{320} holds.\end{pf}

     The fibers of the bundle  $\ADP$  possess  a  natural  quantum
group structure. Further, the bundle $\ADP$ acts on  the  bundle  $P$,
preserving fibers and the right  action.  This  is  a  geometrical
background for the next proposition.

\begin{pro}\label{pro:33} The following identities hold
\begin{gather}
(\id\otimes\phi_M)\phi_M=(\phi_M\otimes \id)\phi_M\\
(\id\otimes F)\beta_M=(\beta_M\otimes \id)F\label{Fcov}\\
(\id\otimes\beta_M)\beta_M=(\phi_M\otimes \id)\beta_M\label{322}\\
(\id\otimes \e_M)=(\e_M\otimes \id)\phi_M=\id \label{324}  \\
(\e_M\otimes \id)\beta_M=\id \label{325}     \\
m_M(\k_M\otimes \id)\phi_M=m_M(\id\otimes \k_M)\phi_M=j_M\e_M \label{326}
\end{gather}
where $m_M \colon \cal{D}\otimes_M \cal{D}\rightarrow  \cal{D}$  is
the multiplication map.
\end{pro}
\begin{pf} In terms of local trivializations, everything  reduces
to elementary algebraic properties of  the  comultiplication,  the
counit, and the antipode. \end{pf}

     We pass to the analysis  of  gauge  transformation,  in  this
quantum framework. In analogy with classical  geometry,  these
transformations will be defined as vertical automorphisms  of  the
bundle.

 \begin{defn}\label{defn:32}
 {\it A gauge transformation} of the  bundle  $P$
 is  every
$S(M)$-linear *-automorphism
$\gamma\colon \cal{B}\rightarrow \cal{B}$ such that
the diagram
\begin{equation}\label{d1}
\begin{CD}
\cal{B} @>{\mbox{$F$}}>> \cal{B}\otimes\cal{A}\\
@V{\mbox{$\gamma$}}VV @VV{\mbox{$\gamma\otimes \id$}}V\\
\cal{B} @>>{\mbox{$F$}}> \cal{B}\otimes\cal{A}
\end{CD}
\end{equation}
is commutative.
\end{defn}
The above diagram says that $\gamma$ intertwines  the  right
action of $G$ on $P$, while the $S(M)$-linearity property ensures
that $\gamma$ is a ``vertical'' automorphism of $P$.
Obviously, gauge transformations form a subgroup
$\cal{G}\subseteq \mbox{Aut}(\cal{B})$.

\begin{pro} (i) The formula
\begin{equation}\label{327}
     f \leftrightarrow  (f\otimes \id)\beta_M =\gamma
\end{equation}
establishes a bijection  between
$S(M)$-linear *-homomorphimsms
$f\colon \cal{D}\rightarrow S(M)$
and gauge transformations $\gamma\in \cal{G}$.
In  terms  of this correspondence, the map $\e_M$
corresponds to the neutral element
in $\cal{G}$ while the product and the inverse  in  the  gauge  group  are
given by
\begin{gather}
            f\k_M  \leftrightarrow  \gamma^{-1} \label{329}\\
  (f'\otimes f)\phi_M  \leftrightarrow  \gamma\gamma'. \label{330}
  \end{gather}

(ii) Let $\gamma$ be an arbitrary gauge  transformation.  Then  the  map
$\gamma_{cl}\colon P_{cl}  \rightarrow P_{cl}$   defined by
\begin{equation} \label{331}
\gamma_{cl}  (p)=p\gamma^{-1}
\end{equation}
is an ordinary gauge transformation of $P_{cl}$.  Moreover,  the  above
formula  establishes  an  isomorphism  between  groups  of   gauge
transformations of bundles $P$ and $P_{cl}  . $ \end{pro}

 \begin{pf}  Identity \eqref{Fcov}  implies  that  a $S(M)$-linear
homomorphism
$\gamma\colon
\cal{B}\rightarrow \cal{B}$ given by the right-hand side of \eqref{327}
satisfies \eqref{d1}.  Identity
\eqref{325} ensures that $\e_M$  corresponds to the neutral element of
$\cal{G}$.

Let us consider an arbitrary  gauge  transformation
$\gamma\in  \cal{G}$.  In
terms of the trivialization system $\tau$ we have
     \begin{equation}\label{332}
     \pi_U \gamma(b)=\sum_{i} \varphi_i \gamma_U (a^{(1)}_i )\otimes
a^{(2)}_i     \end{equation}
for each $U\in \cal{U}$. Here, $\pi_U (b)=\Sum_{i} \varphi_i
\otimes  a_i $  and
$\gamma_U \colon U\rightarrow G_{cl}$   are smooth functions uniquely
determined
by $\gamma$ (understood here in the ``dual'' manner). We have
\begin{equation}\label{333}
\bigl(\gamma_V (a^{(1)}){\restr}_{U\cap V}\bigr)
g_{VU}(a^{(2)})=\gamma_U(a){\restr}_{U\cap V}
       \end{equation}
for each $a\in \cal{A}$ and  $(U,V)\in  N^2 (\cal{U})$.

Conversely,  if
*-homomorphisms $\gamma_U \colon  \cal{A}\rightarrow S(U)$
are  given  such   that
equalities \eqref{333} hold then formula \eqref{332} consistently determines
a gauge transformation $\gamma$.

 Let   us   now    consider    a    map    $f\colon    \Sigma
(\cal{U})\rightarrow S(\cal{U})$ defined by
$$  f=\sideset{}{^\oplus}\sum_{U\in\cal{U}}
                         {f}_U ,$$
where ${f}_U \colon S(U)\otimes \cal{A}\rightarrow S(U)$
are maps  given
by ${f}_U (\varphi\otimes a)=\varphi \gamma_U (a)$.
It is easy to  see  that  if
$b\in \cal{D}$ then $f(b)\in S(M)$
(where  $S(M)$  is  understood  as  a
subalgebra of   $S(\cal{U})$).
Let   us pass to the  corresponding  restriction $f\colon
\cal{D}\rightarrow  S(M)$.  By
construction  \eqref{327}   holds   (it   is   evident   in   a   local
trivialization). Conversely, if  $f\colon  \cal{D}\rightarrow  S(M) $
determines      a      gauge      transformation      $\gamma$       then
$$f(b){\restr}_U =\sum_{i}\varphi_i \gamma_U(a_i ).$$
This easily follows from \eqref{327}.

Let us check correspondences \eqref{329}--\eqref{330}. We have
   \begin{align*}
   &\left[(f\otimes f')\phi_M \otimes \id\right]\beta_M
   =(f\otimes \gamma')\beta_M =\gamma'\gamma,  \\
    & (f\k_M \otimes f)\phi_M =fm_M (\k_M \otimes \id)\phi_M =\e_M .
    \end{align*}
Finally,  the  second  statement  easily  follows  from   the
definition of gauge transformations, and from the local expression
\eqref{332} for them. \end{pf}

     A geometrical explanation of  the  statement  ({\it ii\/})  is  this.
Gauge  transformations,  being  diffeomorphisms  of   $P$   at   the
geometrical level, must preserve classical and quantum parts of $P$.
On the other hand, because of  the  intertwining  property,  gauge
transformations   $\gamma$   are   completely   determined    by    their
``restrictions'' $\gamma_{cl}$ on $P_{cl}$,
which correspond precisely to  the standard
gauge transformations of $P_{cl}$.

 The  quantum  gauge  bundle   $\ADP $  is   also   an   inherently
inhomogeneous geometrical object. This is  a  consequence  of  the
inhomogeneity  of  $G$.  The  classical  part  of  the  bundle  $\ADP$
(*-characters on $\cal{D}$)  is  naturally  identificable  with  the
ordinary gauge bundle of $P_{cl}$. In other words,
$$(\ADP)_{cl}=\AD(P_{cl}  ).$$
Let  $f\colon  \cal{D}\rightarrow  S(M) $  be   the   *-epimorphism
corresponding to $\gamma\in \cal{G}$. This map determines a section $f^*$
of  the
bundle $(\ADP)_{cl}$ as follows
\begin{equation}\label{335}
\left[f^* (x)\right](\varphi)=\left[f(\varphi)\right](x)
\end{equation}
where $x\in M$ and $\varphi\in  \cal{D}$.  In  the  framework  of  the
correspondence  \eqref{327},   the   map  $ f^*$    becomes   the   section
corresponding  to   the   gauge
transformation $\gamma_{cl}$, in   the   classical   manner.

     ****

 We pass  to  the  construction  and  the  study  of  differential
calculus on the bundle $\ADP$. The calculus will be constructed  by
combining differential  forms  on  the  base  manifold  $M$  with  a
differential calculus on the quantum group $G$. This calculus
will be  based  on
the  universal  differential  envelope  $\Gamma^{\wedge}$
of  the   minimal
admissible   first-order   bicovariant *-calculus
$\Gamma$ over $G$.

 \begin{lem} \label{lem:35}  (i)  For  each
 $(U,V)\in  N^2 (\cal{U})$ there exists the   unique   homomorphism
$\xi^{\wedge}_{UV}\colon\Omega(U\cap V)\grten
\Gamma^{\wedge} \rightarrow    \Omega(U\cap V)\grten
\Gamma^{\wedge}$  of  (graded)  differential  algebras,
extending  the  map
$\xi_{UV}$. The map $\xi^{\wedge}_{UV}$
is *-preserving and bijective, and
\begin{equation}\label{336}
(\xi^{\wedge}_{UV})^{-1}=\xi^{\wedge}_{VU}.
\end{equation}

 (ii) We have
\begin{equation}\label{337}
   \xi^{\wedge}_{UV}    \xi^{\wedge}_{VW}(\varphi)
=\xi^{\wedge}_{UW}(\varphi)
\end{equation}
for  each  $(U,V,W)\in N^3(\cal{U})$ and $\varphi\in\WC(U{\cap}V{\cap} W)
\grten\Gamma^\wedge$.

   (iii) The diagrams
\begin{gather*}
\begin{CD}
\Omega(U\cap V)\grten\Gamma^\wedge @>{\mbox{$\id\otimes
\widehat{\phi}$}}>>
\bigl[\Omega(U\cap V)\otimes\Gamma^\wedge
\bigr]\grten_{U\cap V}\bigl[\Omega(U\cap V)\otimes\Gamma^\wedge\bigr]\\
@V{\mbox{$\xi_{UV}^\wedge$}}VV @VV{\mbox{$\xi_{UV}^\wedge
\otimes
\xi_{UV}^\wedge $}}V\\
\Omega(U\cap V)\grten\Gamma^\wedge @>>{\mbox{$\id\otimes
\widehat{\phi}$}}>
\bigl[\Omega(U\cap V)\otimes\Gamma^\wedge
\bigr]\grten_{U\cap V}\bigl[\Omega(U\cap V)\otimes\Gamma^\wedge\bigr]
\end{CD}\\
\begin{CD}
\Omega(U\cap V)\grten\Gamma^\wedge
@>{\mbox{$\id\otimes \e\Pi$}}>> \Omega(U\cap V)\\
@V{\mbox{$\xi_{UV}^\wedge$}}VV @VV{\mbox{$\id$}}V\\
\Omega(U\cap V)\grten\Gamma^\wedge @>>{\mbox{$\id\otimes \e\Pi$}}>
\Omega(U\cap V)
\end{CD}\qquad
\begin{CD}
\Omega(U\cap V)\grten\Gamma^\wedge
@>{\mbox{$\id\otimes\widehat{\kappa}$}}>>
\Omega(U\cap V)\grten\Gamma^\wedge\\
@V{\mbox{$\xi_{UV}^\wedge$}}VV @VV{\mbox{$\xi_{UV}^\wedge$}}V\\
\Omega(U\cap V)\grten\Gamma^\wedge
@>>{\mbox{$\id\otimes\widehat{\kappa}$}}>
\Omega(U\cap V)\grten\Gamma^\wedge
\end{CD}\\
\begin{CD}
\Omega(U\cap V)\grten\Gamma^\wedge @>{\mbox{$\id\otimes
\widehat{\phi}$}}>>
\bigl[\Omega(U\cap V)\otimes\Gamma^\wedge
\bigr]\grten_{U\cap V}\bigl[\Omega(U\cap V)\otimes\Gamma^\wedge\bigr]\\
@V{\mbox{$\psi_{UV}^\wedge$}}VV @VV{\mbox{$\xi_{UV}^\wedge
\otimes
\psi_{UV}^\wedge $}}V\\
\Omega(U\cap V)\grten\Gamma^\wedge @>>{\mbox{$\id\otimes
\widehat{\phi}$}}>
\bigl[\Omega(U\cap V)\otimes\Gamma^\wedge
\bigr]\grten_{U\cap V}\bigl[\Omega(U\cap V)\otimes\Gamma^\wedge\bigr]
\end{CD}
\end{gather*}
are commutative.
\end{lem}

\begin{pf} The uniqueness of $\xi^{\wedge}_{UV}$
follows  from  the  fact
that $\WC(U\cap V)\grten \Gamma^{\wedge}$
is  generated,  as  a  differential
algebra, by $\SC (U\cap V)\otimes \cal{A}$.  The hermicity  of
$\xi^{\wedge}_{UV}$      follows
from the fact that $*\xi^{\wedge}_{UV} *$
is a differential  extension  of  the
same map $*\xi_{UV} *=\xi_{UV}$. In a similar way it
follows from Lemma~\ref{lem:31}
that above diagrams
are commutative, and  that
\eqref{336}--\eqref{337} hold.

 We prove the existence of $\xi^{\wedge}_{UV}$.
The admissibility of $\Gamma$ and
the universality of $\Gamma^{\wedge}$
imply that maps $g_{UV}$ admit the unique
graded-differential    (*-preserving)     extensions $\widehat{g}_{UV}
\colon\Gamma^{\wedge} \rightarrow  \Omega(U\cap  V)$.

Now,  the  maps   $f_{UV}  \colon
\Gamma^{\wedge} \rightarrow \Omega(U\cap V)\grten\Gamma^{\wedge}  $
given by
$$ f_{UV}(w)=\sum_i \widehat{g}_{UV}
(\widehat{\k}(w^1_i )w^3_i )\otimes w^2_i ,$$
where
$\Sum_i w^1_i \otimes w^2_i \otimes w^3_i
=(\widehat{\phi}\otimes \id)\widehat{\phi}(w)
=(\id\otimes \widehat{\phi})\widehat{\phi}(w
),$ are  homomorphisms  of  differential  *-algebras.  Finally  let
$\xi^{\wedge}_{UV}$     be defined by
$$  \xi^{\wedge}_{UV}(\alpha\otimes w)=\alpha f_{UV}  (w).$$
It is evident that such  defined  maps  are  differential  algebra
homomoprhisms extending $\xi_{UV}$.\end{pf}

 Let $\Omega(\tau,\ADP)$ be the set of  all  elements
 $w\in\Sigma^{\wedge}
 (\cal{U})$ satisfying
 \begin{equation}  \label{342}
 (_U {\restr}_{U\cap V}   \otimes \id)p_U (w)=\xi^{\wedge}_{UV}
 (_V {\restr}_{U\cap V}   \otimes \id)p_V (w),    \end{equation}
for each $(U,V)\in N^2 (\cal{U})$. It is clear that $\Omega(\tau,\ADP) $
is a graded- differential *-subalgebra of $\Sigma^{\wedge}   (\cal{U}),$
and
that $\Omega^0 (\tau,\ADP)=\cal{D}$. The  elements  of  the  algebra
$\Omega(\tau,\ADP)$ play the role  of  differential  forms  on  the
bundle $\ADP$. This algebra is generated by $\cal{D},$  and  in  fact
does not depend of a trivialization system $\tau$.  More  precisely,
if $\eta$ is another trivialization system for $P$ then  there  exists
(the       unique)       differential      $ (*-)$       isomorphism
$\Omega(\tau,\ADP)\leftrightarrow  \Omega(\eta,\ADP)$  extending  the
identity map on $\cal{D}$. For this reason  we  shall  simply  write
$\Omega(\tau,\ADP)=\Omega(\ADP)$.

\begin{pro}\label{pro:36}
(i) The maps  $\bigl\{\e_M, j_M, \beta_M, \phi_M\bigr\}$   admit  unique
extensions
\begin{gather*}
\phi^{\wedge}_M \colon\Omega(\ADP)\rightarrow \Omega(\ADP)
\grten_M  \Omega(\ADP)\\
\e^{\wedge}_M \colon \Omega(\ADP)\rightarrow \Omega(M)\\
j^{\wedge}_M \colon \Omega(M)\rightarrow \Omega(\ADP)\\
\beta^{\wedge}_M \colon \Omega(P)\rightarrow \Omega(\ADP)
\grten_M  \Omega(P),
\end{gather*}
which are homomorphisms of graded-differential algebras.

(ii) The map $\k_M$ admits the unique extension
$\k^{\wedge}_M \colon \Omega(\ADP)\rightarrow \Omega(\ADP)$ which is
graded-antimultiplicative and satisfies
\begin{equation}\label{343}
\k^{\wedge}_M d=d\k^{\wedge}_M .
\end{equation}

   (iii) The following identities hold
   \begin{gather}
 (\phi^{\wedge}_M \otimes \id)\phi^{\wedge}_M
 =(\id\otimes \phi^{\wedge}_M )\phi^{\wedge}_M  \label{344}\\
 (\phi^{\wedge}_M \otimes \id)\beta^{\wedge}_M
 =(\id\otimes \beta^{\wedge}_M )\beta^{\wedge}_M  \label{345}\\
                         (\id\otimes \widehat{F})\beta^{\wedge}_M
                         =(\beta^{\wedge}_M \otimes \id)\widehat{F}
\label{346}\\
 (\id\otimes
\e^{\wedge}_M )\phi^{\wedge}_M =(\e^{\wedge}_M \otimes
\id)\phi^{\wedge}_M =\id \label{347}\\
                          (\e^{\wedge}_M \otimes \id)\beta_M =\id
\label{348}\\
      m^{\wedge}_M (\k^{\wedge}_M \otimes \id)\phi^{\wedge}_M
      =m^{\wedge}_M (\id\otimes \k^{\wedge}_M )\phi^{\wedge}_M
      =j^{\wedge}_M \e^{\wedge}_M , \label{349}\end{gather}
where $m^{\wedge}_M$  is the multiplication map in $\Omega(\ADP). $

(iv) We have
\begin{equation}\label{350}
       *\k^{\wedge}_M *=(\k^{\wedge}_M )^{-1},
\end{equation}
while $\bigl\{\e^{\wedge}_M,j^{\wedge}_M,
\beta^{\wedge}_M,\phi^{\wedge}_M\bigr\}$ are *-preserving maps.
\end{pro}
 \begin{pf}  Using  the   (anti)multiplicativity,   the   intertwining
differentials  properties,  and  the  fact  that  all   considered
differential algebras are generated by corresponding zero-th order
subalgebras, it is easy to see that extensions of all maps in  the
game are, if exist, unique. The  same  properties,  together  with
Proposition~\ref{pro:33},  imply  that  identities
\eqref{344}--\eqref{349}  hold.   The
statement ({\it iv\/}) follows from
({\it ii\/}) Proposition~\ref{pro:32}
in a similar way.
Finally,  existence  of  maps  $\e^{\wedge}_M$,
$j^{\wedge}_M$, $\beta^{\wedge}_M$, $\phi^{\wedge}_M$
and   $\k^{\wedge}_M$    can   be
established in a similar way as for maps $\e_M$, $j_M$,
$\beta_M$, $\phi_M$  and $\k_M$. \end{pf}

Every $\gamma\in \cal{G}$ understood as a  *-homomorphism
$f\colon  \cal{D}\rightarrow  S(M)$  is  uniquely  extendible  to  a
$\Omega(M)$-linear *-homomorphism
$f^{\wedge} \colon  \Omega(\ADP)\rightarrow  \Omega(M)$
of  graded-differential  algebras.

The following correspondences hold
 \begin{align}
 \gamma^{-1}&   \leftrightarrow f^{\wedge} \k^{\wedge}_M
 \label{351}\\
 \gamma\gamma' &\leftrightarrow m^{\wedge}_M (f^{\wedge}
\otimes {f}^{\prime\wedge} )\phi^{\wedge}_M.  \label{352}\end{align}

\section{Gauge Fields}
In this section we shall present a generalization of the
classical gauge theory, within the geometrical framework of
quantum principal bundles. The base manifold $M$ will play the role
of space-time. The quantum group $G$ will describe ``internal
symmetries'' of the system. In order to simplify considerations,
we shall deal only with a ``pure gauge theory''.

Let us assume that $\Gamma_{\inv}$ is endowed with an $\Ad$-invariant
scalar product $(,)$. It means that

$$(\vartheta ,\eta )\otimes 1=\sum_{kl} (\vartheta_k ,\eta_l )
\otimes c^*_k d_l  $$
for each $\vartheta ,\eta \in\Gamma_{\inv} ,$ where
$\Sum_k \vartheta_k \otimes c_k =\Ad(\vartheta )$
and $\Sum_l\eta_l\otimes d_l =\Ad(\eta )$.

Let us assume that M is oriented and endowed with a
(pseudo)riemannian structure.

Let us denote by ${\star}$ the Hodge operation on $\Omega (M)$.
It can be (uniquely) extended to a linear map ${\star}\colon
\hor(P)\rightarrow\hor(P)$ such that

$${\star}\bigl(i^{\wedge} (\alpha )b\bigr)=i^{\wedge} \bigl({\star}(\alpha )
\bigr)b, $$
for each $\alpha\in \Omega (M)$ and $b\in \cal{B}$.

Following the classical analogy gauge fields will be
geometrically represented by connection forms $\omega$ on the bundle $P$.

To make possible dynamical considerations it is necessary to
fix a lagrangian. Generalizing the classical situation, it is
natural to consider lagrangians which are quadratic functions of
the curvature $R_{\omega} $. The curvature operator $R_{\omega}$
depends, besides on
the connection $\omega ,$ also on a choice of the embedded differential
map $\delta \colon\Gamma_{\inv} \rightarrow\Gamma_{\inv}
\otimes\Gamma_{\inv}$. As a consequence of this, dynamical properties of the
gauge theory
will be essentially influenced by $\delta$. In the classical case the
curvature is $\delta$-independent.

Let us consider a map $L\colon\con(P)\rightarrow\hor(P)$ given by
\begin{equation}\label{41}
L(\omega )=\sum_iR_{\omega}(e_i){\star}\left[R_{\omega}(\bar{e}_i)\right],
\end{equation}
where elements $e_i$ form an orthonormal system in $\Gamma_{\inv}$ and the bar
denotes the conjugation in $\Gamma_{\inv}$.
It is easy
to see that $L(\omega )$ is independent of the choice of the mentioned
orthonormal system.

The map $L$ in fact takes values from the space $\Omega^n (M)$ (where
$n$ is the dimension of $M$). Indeed, in terms of local
trivializations we have
\begin{equation}\label{42}
\pi^{\wedge}_U \bigl[L(\omega )\bigr]=\sum_iF^U (e_i ){\star}
\bigl[F^U (\bar{e}_i )\bigr]\otimes 1. \end{equation}
This easily
follows from the fact that $\Sum_i e_i\otimes
\bar{e}_i$ is $\Ad^{\otimes 2}$-invariant.

We shall interpret the map $L$ as the lagrangian. In terms of the local
representation, stacionary points of the corresponding action functional
$S(\omega )=\int_M L(\omega )$ are given by the following equations of
motion
\begin{equation}\label{45}
d{\star} F^U (\bar{e}_k )+\frac{1}{2}\sum_{ij}(d^{jk}_i -d^{kj}_i )A^U(e_j)
{\star} F^U(\bar{e}_i)=0
\end{equation}
where numbers $d^{ij}_k$ are determined by
\begin{equation}\label{44}
\delta (e_k )=-\frac{1}{2} \sum_{ij}d^{ij}_k e_i \otimes e_j.
\end{equation}

The above equations correspond to the classical Yang-Mills
equations of motion. The numbers $(d^{jk}_i -d^{kj}_i)/2$
play the role of the
structure constants of (the Lie algebra of) $G$.

If the space $\Gamma_{\inv}$ is infinite-dimensional a technical
difficulty arises, related to a question of convergence of the sum
in \eqref{41}--\eqref{42}. In such cases, it is necessary
to restrict possible
values of $\omega$ on some subspace of $\con(P),$ consisting of connections
having sufficiently rapidly decreasing components, in an
appropriate sense.

****

We pass to the study of symmetry properties of the introduced
lagrangian.
As first, it is easy to see that $L(\omega )$ is invariant under
gauge transformations of the bundle $P. $

The group $\cal{G}$ naturally acts on the left, via compositions,
on the space $\psi (P)$ of
pseudotensorial forms.
The space $\tau (P)$ is invariant with respect to this action, because
$\hor(P)$ is $\cal{G}$-invariant.
The connection space is gauge invariant, too.
In terms of gauge potentials the transformation of connections is
\begin{equation}\label{46}
A^U (\vartheta )\longrightarrow \sum_kA^U (\vartheta_k )\gamma^U(c_k)
+\partial^U(\vartheta ).
\end{equation}
Here, $\Ad(\vartheta )=\Sum_k\vartheta_k \otimes c_k$,
the map $\partial^U\colon\Gamma_{\inv} \rightarrow\Omega^1 (U)$ is
given by
$$\partial^U\pi (a)=\gamma^U\!\k(a^{(1)})d\gamma^U(a^{(2)}),$$
while $\gamma^U\colon \cal{A}\rightarrow S(U)$ is the map locally
representing $\gamma$. Further, the
transformation of the curvature is
\begin{equation}\label{47}
F^U (\vartheta )\longrightarrow
\sum_kF^U (\vartheta_k )\gamma^U(c_k).   \end{equation}

The lagrangian \eqref{41} is invariant under gauge transformations
of the bundle $P$. This is a simple consequence of the unitarity of
the representation $\Ad$.

This invariance is a manifestation of {\it classical symmetry
properties} of the lagrangian. These symmetry properties are
completely expressible in terms of the classical part $P_{cl}$ of $P$.

On the other hand, the lagrangian $L(\omega )$ possesses symmetry
properties which are not expressible in classical terms. The
appearance of these ``quantum symmetries''  is a purely quantum
phenomena caused by the quantum nature of the space $G$. Formally,
they can be described as the invariance of the lagrangian under a
natural action of the quantum gauge bundle $\ADP$.

Let $\psi (P,\ADP)$ be the space of linear maps
$f\colon\Gamma_{\inv} \rightarrow\Omega (\ADP)\grten_M
\Omega (P)$ satisfying
\begin{equation}\label{48}
(f\otimes\id)\varpi=(\id\otimes F^\wedge)f.
\end{equation}
If $\varphi \in\psi (P)$ then $\beta^{\wedge}_M
\varphi \in\psi (P,\ADP)$. Hence
it is possible
to introduce the map $\beta^{\wedge}_M \colon\psi (P)
\rightarrow\psi (P,\ADP)$ (via compositions).

Let us compute the element $\beta^{\wedge}_M \omega$ for
$\omega \in \con(P)$. Using the
definition of $\beta^{\wedge}_M$ and the local expression for
$\omega$ we obtain
\begin{multline}\label{49}
(\pi_U^\sharp\otimes\pi_U)\bigl[\beta^{\wedge}_M (\omega )(\vartheta )\bigr]
=1_U \otimes 1\otimes\vartheta\\
{}+\sum_k\Bigl\{A^U (\vartheta
_k)\otimes c^{(1)}_k\otimes c^{(2)}_k +1_U \otimes\vartheta_k
\otimes c_k \Bigr\}. \end{multline}

Here an identification
$$\bigl[\Omega (U)\grten\Gamma^{\wedge} \bigr]
\grten_U\bigl[\Omega (U)\grten\Gamma^{\wedge}\bigr]
=\Omega (U)\grten\Gamma^{\wedge}
\grten\Gamma^{\wedge}$$
is assumed.
It is worth noticing that the transformation law \eqref{46} is
contained in \eqref{49}. Indeed, understanding gauge transformations as
differential algebra homomorphisms $f^\wedge\colon\Omega(\ADP)
\rightarrow\Omega (M)$ we obtain \eqref{46}
by composing $\beta^{\wedge}_M (\omega )$ and $f^\wedge\otimes \id$.

The curvature is transformed as follows
\begin{equation}\label{410}
(\pi_U^\sharp\otimes\pi_U)\bigl[\beta^{\wedge}_M (R_{\omega} )
(\vartheta )\bigr]=\sum_kF^U (\vartheta_k)\otimes c^{(1)}_k \otimes
c^{(2)}_k .  \end{equation}
That is, in local terms we have
\begin{equation}\label{411} F^U \longrightarrow (F^U \otimes \id)\Ad.
\end{equation}
The curvature operator is gauge-covariant in the sense that
\begin{equation}\label{412}
\beta^{\wedge}_M (R_{\omega})=d\beta^{\wedge}_M (\omega )
-\langle \beta^{\wedge}_M(\omega),\beta^{\wedge}_M(\omega)\rangle.
\end{equation}

A possible interpretation of the above equation (which is a
trivial consequence of the fact that $\beta^{\wedge}_M \colon\Omega (P)
\rightarrow\Omega (\ADP)\grten_M\Omega
(M)$ is a
differential algebra homomorphism) is this. The relation between
the connection $\omega$ and its curvature is, being expressible in
intrinsicly geometrical terms, preserved under the action of $\ADP$.
Expression \eqref{47} for the curvature of the transformed connection
under a gauge transformation also directly follows from \eqref{410}.

In order to find the transformation of the local expression
for the lagrangian, we should insert in \eqref{42} the local
expression for the transformed curvature, under the action
$\beta^{\wedge}_M$
of $\ADP$ on $P$. The lagrangian transforms as follows
\begin{equation}\label{413}
\Bigl\{\sum_kF^U (e_k){\star} F^U(\bar{e}_k)\Bigr\} \longrightarrow
\sum_{kln}F^U (e_l){\star} F^U (\bar{e}_n)\otimes c_{lk} c^*_{nk}
\end{equation}
where $\ad(e_i)=\Sum_j e_j\otimes c_{ji}$. On the other hand
\begin{equation}\label{414}
\Bigl[\sum_kF^U(e_k){\star} F^U(\bar{e}_k)\Bigr]\otimes 1=
\sum_{kln}F^U (e_l){\star} F^U (\bar{e}_n)\otimes c_{lk} c_{nk}^*,
\end{equation}
because of the $\Ad^{\otimes 2}$-invariance of
$\Sum_ke_k\otimes \bar{e}_k$.

Hence, the lagrangian is invariant with respect to the action
$\beta^{\wedge}_M$ of the gauge bundle $\ADP$ on $P$.

It is important to mention that the property of ``quantum
gauge invariance'' of the lagrangian can not be viewed as an
inherent property of the local expression \eqref{42}. Because this
property essentially depends on {\it the ordering of terms}
$F^U$ and ${\star} F^U$.
However, in the general case, the ordering of terms $R_{\omega}$
and ${\star} R_{\omega}$ in
the global representation of the lagrangian is essential, because
$\cal{B}$ is a noncommutative algebra.
\section{Example}
We shall now illustrate the presented formalism on a concrete
example, assuming that $G=SU_{\mu} (2)$ (with $\mu =(-1,1)\setminus\{0\}$).
By definition
\cite{W1} this compact matrix quantum group is based
on the $2\times 2$ matrix
\begin{equation}\label{51}
u=\begin{pmatrix}\alpha & -\mu\gamma^*\\
\gamma &\phantom{-\mu}\alpha^* \end{pmatrix}\end{equation}
where the elements $\alpha$ and $\gamma$ satisfy the following relations
\begin{equation}\label{52}\begin{gathered}
\alpha \gamma =\mu \gamma \alpha \quad\gamma\alpha^*=\mu\alpha^*\gamma
\quad
\gamma \gamma ^*=\gamma^*\gamma \\
\alpha^*\alpha +\gamma^*\gamma =1 \qquad\alpha \alpha^*+\mu^2
\gamma^* \gamma =1. \end{gathered}\end{equation}

The classical part of $G$ is isomorphic to $U(1)$. An explicit
isomorphism is given by $g\leftrightarrow g(\alpha)$.

It turns out (cite{D}--Section 6) that the right $\cal{A}$-ideal
$\widehat{\cal{R}}\subseteq \ker(\e)$ determining the
minimal admissible (bicovariant $*-$) first-order calculus $\Gamma$
over $G$ is given by
\begin{equation}\label{53}
\widehat{\cal{R}}=\bigl(\mu^2\alpha +\alpha^*-(1+\mu^2)1\bigr)\ker(\e).
\end{equation}
Let $X\colon \cal{A}\rightarrow \Bbb{C}$ be a generator of
$\lie(G_{cl} )$ specified by
\begin{equation}\label{54}\begin{gathered}
X(\alpha )=-X(\alpha^*)=1/2\\
X(\gamma )=X(\gamma^*)=0.\end{gathered}\end{equation}

Let $\rho\colon \cal{A}\rightarrow\cal{A}$ be a map given by
$\rho =(X\otimes \id)\ad$. Let $\nu \colon\Gamma_{\inv}\rightarrow \Bbb{C}$ and
$\widetilde{\rho}\colon \Gamma_{\inv}\rightarrow\cal{A}$
be the maps defined by $\nu \pi =X$ and $\widetilde{\rho} \pi =\rho $. Then
$\widetilde{\rho} =(\nu\otimes \id)\Ad,$
and $\widetilde{\rho}$ maps isomorphically the space
$\Gamma_{\inv}$ onto the *-subalgebra $\cal{Q}\subseteq\cal{A}$
of left $G_{cl}$-invariant elements of $\cal{A}$. The subalgebra $\cal{Q}$ is
interpretable as the algebra of polynomial functions on a
quantum $2$-sphere.

The adjoint action $\Ad$ is reducible. The space $\Gamma_{\inv}$ is
decomposable into the orthogonal sum
$$\Gamma_{\inv}=\sideset{}{^\oplus}\sum_{k\geq 0}\Gamma^k_{\inv}$$
of irreducible
subspaces. The subspace $\Gamma^k_{\inv}$ is $(2k+1)$-dimensional
(that is, all integer-spin irreducible multiplets are
in the game).

The space $\widetilde{\rho} (\Gamma^k_{\inv})=\cal{Q}_k$
is spanned by quantum spherical harmonics
$\zeta_{km}$, where $m\in\{-k,\dots, k\}$. They
constitute a standard basis for the
action of $G$. Explicitly, these elements are given by
\begin{equation}\label{55}\begin{aligned}
&\zeta_{km}=(-)^m \mu^{km-m}\bigl[(k-m)_{\mu}!/(k+m)_{\mu}!\bigr]^{1/2}
\partial^mp_k (\gamma \gamma^*)\gamma^m\alpha^m \\
&\zeta_{k,-m}= \mu^{km}\alpha^{*m}\gamma^{*m}\bigl[(k-m)_{\mu}!/(k+m)_{\mu}
!\bigr]^{1/2}\partial^mp_k (\gamma
\gamma^*).\end{aligned}\end{equation}

Here, $m\geq 0$ and $\partial\colon P(x)\rightarrow P(x)$
is a ``quantum differential''  (acting on
the space $P(x)$ of $x$-polynoms) specified by $\partial (x^n)=n_{\mu}
x^{n-1} $.
Finally, $p_k (x)$ are polynomials given by
\begin{equation}\label{56}\begin{aligned}
&p_k (x)=(-)^k c_k \partial^k\Bigl[x^k\prod^k_{j=1}(1-\mu^{1-j}x)\Bigr]\\
&p_0 (x)=1, \end{aligned}\end{equation}
while $c_k>0$ and $$k_{\mu}!=\prod^k_{j=1}j_{\mu} \quad\quad
j_{\mu} =\frac{1-\mu^{2j}}{1-\mu^{2\phantom{j}}}.$$

Let us now describe a construction of the natural embedded
differential map $\delta $. We shall first construct a complement
$\cal{L}\subseteq \ker(\e)$
of the space $\widehat{\cal{R}}$.

The elements $\gamma^k$ ($k\in\Bbb{N} )$ are primitive for
the adjoint action of
$G$ on $\ker(\e)$. Let $\cal{L}\subseteq \ker(\e)$
be the minimal $\Ad$-invariant subspace
containing these elements, and the \ad-invariant element
$\mu^2 \alpha +\alpha^* -(1+\mu^2)1$. It turns out that the restriction
$(\pi{\restr}\cal{L})\colon\cal{L}\rightarrow
\Gamma_{\inv}$ is
bijective. Evidently, this restriction intertwines the adjoint
actions. Let $\delta \colon\Gamma_{\inv} \rightarrow\Gamma_{\inv}
\otimes\Gamma_{\inv}$ be defined by
\begin{equation}\label{57}
\delta (\vartheta )=-(\pi \otimes\pi )\phi\left[(\pi{\restr}\cal{L})^{-1}
(\vartheta)\right].
\end{equation}

It is clear, by construction, that $\delta$ is an embedded
differential map. Moreover,
\begin{equation}\label{58}
\delta \varkappa =-(\varkappa \otimes\varkappa )\delta \end{equation}
where the extension of the antipode
$\varkappa \colon\Gamma_{\inv}\rightarrow \Gamma_{\inv}$ is
given by
\begin{equation}\label{59} \varkappa \pi =-\pi \k^2. \end{equation}

Let us compute the values of $\delta$ on the singlet and the triplet
subspace of $\Gamma_{\inv}$. The singlet space $\Gamma^0_{\inv}$ is spanned by
\begin{equation}\label{510}
\varepsilon =\pi (\mu^2 \alpha +\alpha^* ), \end{equation}
while the triplet space $\Gamma^1_{\inv}$ is spanned by
\begin{equation}\label{511}
\eta_+=\pi(\gamma),\quad \eta =\pi (\alpha -\alpha^* ),
\quad \eta_-=\pi (\gamma^* ). \end{equation}

Applying the definition of $\delta$ we obtain
\begin{gather*}
-\delta (\varepsilon )=\bigl(\varepsilon \otimes\varepsilon +\mu^2
\eta \otimes\eta \bigr)/(1+\mu^2 )-
\mu \eta_+\otimes \eta_- -\mu^3 \eta_- \otimes\eta_+ \\
-\delta (\eta_+ )=\bigl((\varepsilon -\mu^2 \eta )\otimes\eta_+
+\eta_+ \otimes(\varepsilon +\eta )
\bigr)/(1+\mu^2 )\\
-\delta (\eta_-)=\bigl(\eta_-\otimes (\varepsilon -\mu^2 \eta )
+(\varepsilon +\eta )\otimes\eta_-
\bigr)/(1+\mu^2 )\\
-\delta (\eta )=\bigl(\varepsilon \otimes\eta +\eta \otimes\varepsilon
+(1-\mu^2 )
\eta \otimes\eta \bigr)/(1+\mu^2 )+\mu \bigl(\eta_+\otimes \eta_- -\eta_-
\otimes\eta_+ \bigr)
\end{gather*}

****

The corresponding gauge theory based on the bundle $P$,
calculus $\Gamma$ group $G$ and the lagrangian $L(\omega )$ is essentially
different from the classical gauge theory with $G=SU(2). $

As first, gauge fields possess infinitely many
internal degrees of freedom. In the classical limit $\mu\rightarrow 1$ the
restriction $A^U {\restr}\Gamma^1_{\inv}$ on the triplet subspace
can be interpreted as
a classical $SU(2)$ gauge field. Restrictions on other irreducible
subspaces are classically interpretable as additional vector
fields.

According to the general theory, the connection $A^U$ can be
decomposed into ``classical'' and ``purely quantum'' parts
$$A^U=A_{cl}^{U} +A_\perp^U,$$
where $A_{cl}^U{\restr} \ker(\nu)=0$ and $A_{\perp}^U
(\varepsilon )=0$. The map $A_{cl}^U$ can be
interpreted as a connection on the classical $U(1)$-bundle $P_{cl}$. It
is important to point out that the decomposition $\Gamma_{\inv}=\ker(\nu)+
\Bbb{C}\varepsilon$ is
incompatible with the decomposition of $\Gamma_{\inv}$ into
irreducible multiplets.

Let us compute the singlet and the triplet components of the
curvature. Applying the definition of $\delta$ and using the local
expression of the curvature we find
\begin{gather*}
F^U(\varepsilon )=dA^U(\varepsilon )+\mu (1-\mu^2)A^U(\eta_-)A^U(\eta_+)\\
F^U (\eta_+)=dA^U (\eta_+)+A^U (\eta_+ )A^U (\eta ) \\
F^U (\eta_-)=dA^U (\eta_-)+A^U (\eta )A^U (\eta_-)\\
F^U (\eta )=dA^U (\eta )+2\mu A^U (\eta_+ )A^U (\eta_-).
\end{gather*}

In general, components of the restriction $F^U {\restr}
\Gamma^k_{\inv}$ will be
expressible through fields $A^U (\vartheta ),$
where $\vartheta \in\Gamma^l_{\inv}$ and $1\leq
l\leq k$.

Equations of motion are mutually {\it essentially correlated}.
Indeed, the equation describing the propagation of fields $A^U{\restr}
\Gamma^k_{\inv}$
will generally contain terms of the form $A^U (\vartheta ){\star} F^U
(\eta)$, where
$\vartheta,\eta \in\Gamma^{i,j}_{\inv}$ and $\vert i-k\vert \leq j
\leq i+k$. This easily follows from the definition
of $\delta$. It is interesting to observe that nonsinglet components are
not explicitly influenced by the singlet component $A^U(\varepsilon )$. On the
other hand, the singlet propagation is intertwined only with
$A^U (\eta_\pm)$. Explicitly, $$d{\star} F^U(\varepsilon )=0.$$
\section{Concluding Remarks}
In this study we have assumed that the higher-order differential
calculus on the structure group is based on the corresponding universal
envelope. All constructions can be performed also in the case when the
higher-order calculus is described by the corresponding bicovariant
(braided) external algebra \cite{W3}.

The admissibility assumption for $\Gamma$ ensures full local
trivializability of differential structures on $P$ and $\ADP$. From the
``local'' point of view however, the whole formalism works for an
arbitrary bicovariant *-calculus $\Gamma$.

In summary, physical properties of the presented gauge theory are
essentially influenced by two additional structural elements. As first,
it is necessary to fix a bicovariant *-calculus $\Gamma$ over $G$. This
determines kinematical degrees of freedom. Secondly, the curvature is
determined only after fixing an embedded differential map $\delta$. In
such a way the dynamics becomes $\delta$-dependent.

The same geometrical framework of quantum principal bundles contains
logically inequivalent ways of generalizing classical gauge theory. For
example, in the formulation discussed in \cite{BM} the calculus
on $G$ does not figure explicitly, and the curvature is defined in
a different way. In this formulation internal degrees of freedom
and the curvature are determined after fixing a fundamental
representation of the structure quantum group, which a priori excludes
quantum phenomena appearing at the level of differential calculus.

The presented gauge theory admits a natural generalization to the
completely quantum context, in which the base manifold $M$ is a
noncommutative space.

\appendix
\section{The Minimal Admissible Calculus}

In this Appendix some properties of entities associated with
the minimal admissible calculus $\Gamma$ over the quantum $SU(2)$ group
are collected. In particular, we shall analyze in more details the
structure of the space $\cal{L}$ which determines the embedded
differential map.

For each integer $n\geq 1$ let $u_n$ be the $n\times n$ matrix over
$\cal{A},$
corresponding to the irreducible representation \cite{W1,W2}
of $G$, having the
spin $(n-1)/2$ and
acting in $\Bbb{C}^n$.
Let $\cal{A}_n$ be the lineal spanned by matrix elements of $u_n$. We have
$$ \cal{A}=\sideset{}{^\oplus}\sum_{n\geq 1}\cal{A}_n ,$$
according to the representation theory of $G$. The spaces $\cal{A}_n$
are invariant under the adjoint action of $G$. They are mutually orthogonal,
relative to the scalar product induced by the Haar measure $h\colon\cal{A}
\rightarrow\Bbb{C}$.

In subspaces $\cal{A}_n$ the adjoint action decomposes (without
degeneracy) into irreducible multiplets with spins from the set
$\{0,1,\dots,n-1\}. $
\begin{lem}\label{lem:A1}
Let $\xi \in \cal{A}_n$ be a primitive element for the $k$-spin
subrepresentation of $\ad{\restr} \cal{A}_n$. Then,
\begin{equation}\label{A1}
\xi =p_{kn} (\lambda )\gamma^k  \end{equation}
where $\lambda =\mu \alpha +\mu^{-1} \alpha^*$ and $p_{kn}$
is a polynom of degree $n-k-1$
with real coefficients.
\end{lem}
\begin{pf} From the representation theory of $G$ it follows that
$$ \cal{A}_2\cal{A}_n=\cal{A}_n\cal{A}_2=\cal{A}_{n-1}\oplus \cal{A}_{n+1}$$
for each $n\geq 2$. This implies that $\cal{A}_n\setminus\{0\}$
is consisting of certain polynoms of degree $n-1$ (over generators).
Further, polynoms of degree $k\leq n-1$ form the space
$\Sum_i^{*\oplus}\cal{A}_i$, where $i\leq n$. Also, from the reality
of commutation relations \eqref{52} and the orthogonality of spaces $\cal{A}_n$
it follows that we can write
$$ \cal{A}_n=\cal{A}_n^\Re\oplus i\cal{A}_n^\Re$$
where $\cal{A}_n^\Re$ is consisting of polynoms with real coefficients.

On the other hand, every non-zero element of the form \eqref{A1}
is primitive, and generates an
irreducible $k$-spin multiplet relative to the adjoint representation.
Having in mind the form of the decomposition of $\Ad{\restr} \cal{A}_n$
into irreducible multiplets we conclude that \eqref{A1} covers all
primitive elements of the restriction $\Ad{\restr}\cal{A}_n$.
\end{pf}

Let us assume that polynoms $p_{kn}$ are fixed. For fixed $k$,
polynoms $p_{kn}$ are orthogonal,
with respect to the scalar product
given by
\begin{equation}\label{A2}
(p,q)=h\bigl(q(\lambda )\gamma^k \gamma^{*k} p(\lambda )^* \bigr),
\end{equation}
Let $j\colon \cal{A}\rightarrow\cal{A}$ be the
modular automorphism \cite{W2} corresponding to the Haar
measure. This map is characterized by the identity
\begin{equation}\label{A3} h(ba)=h\bigl(j(a)b\bigr).\end{equation}
In the case of the quantum $SU(2)$ group we have
\begin{equation}\label{A4}\begin{gathered}
j(\gamma )=\gamma
\qquad j(\alpha )=\mu^2 \alpha \\
j(\alpha^* )=\mu^{-2}\alpha^*\qquad
j(\gamma^* )=\gamma^* .\end{gathered}\end{equation}

Applying \eqref{A3}--\eqref{A4} we see that the scalar product
defined in \eqref{A2}
can be rewritten in the form
\begin{equation}\label{A5}
(p,q)=h\bigl[p^* (\lambda )q(\lambda )(\gamma \gamma^* )^k \bigr].
\end{equation}

Now, starting from \eqref{A5}, observing that the above scalar
product is invariant under the replacement
$\lambda\rightarrow \alpha +\alpha^* ,$ and using
elementary properties of polynoms it can be shown that all zeroes
of $p_{kn}$ are contained in the interval $[-2,2]$.

We have
\begin{equation}\label{A6}
\cal{L}=\Bbb{C}\bigl(\mu \lambda -(1+\mu^2 )1\bigr)\oplus
\Bigl\{\sideset{}{^\oplus}\sum_{k\geq 1}\cal{L}_k \Bigr\},
\end{equation}
where $\cal{L}_k\subseteq \cal{A}_{k+1}$ is the $k$-spin
irreducible subspace (for the adjoint
action). Let $\cal{L}_*\subseteq\cal{A}$ be the lineal given by
\begin{equation}\label{A7}
\cal{L}_*=\Bbb{C}1\oplus\Bigl\{\sideset{}{^\oplus}\sum_{k\geq 1}
\cal{L}_k\Bigr\}.\end{equation}

Let $P(\lambda )\subseteq\cal{A}$
be the subalgebra generated by $\lambda$.
\begin{lem}\label{lem:A2} (i) We have
\begin{equation}\label{A8}
\pi (\cal{A}_n)=\sideset{}{^\oplus}\sum_{k\leq n-1} \Gamma^k_{\inv}
\end{equation}
for each $n\in\Bbb{N}$.

(ii) The map $\mho\colon P(\lambda )\otimes \cal{L}_*
\rightarrow\cal{A}$ given by
\begin{equation}\label{A9}
\mho\bigl(p(\lambda )\otimes a\bigr)=p(\lambda )a
\end{equation}
is bijective. \end{lem}

\begin{pf} Let us prove that $\mho$ is bijective. As first, let us
observe that the elements from $P(\lambda )$ are $\ad$-invariant.
In particular,
\begin{equation}\label{A10}
\ad\bigl(p(\lambda )a\bigr)=p(\lambda )\ad(a),
\end{equation}
for each $a\in\cal{A}$. According to Lemma~\ref{lem:A1}
all primitive elements for
$\ad\colon\cal{A}\rightarrow\cal{A}\otimes\cal{A}$
are contained in the image of $\mho$. Now \eqref{A10}
implies that $\mho$ is surjective. We prove that $\mho$
is injective. It is
sufficient to check that $\mho {\restr}\bigl( P(\lambda )\otimes
\cal{L}_k\bigr)$
is injective, for each $k\in\Bbb{N}$.
However, it follows again from Lemma~\ref{lem:A1}
and \eqref{A10}, because
$\mho\bigl(p_{kn} (\lambda )\otimes\cal{L}_k\bigr)\subseteq\cal{A}_n$
is exactly the $k$-spin irreducible subspace.

The following identity holds
\begin{equation}\label{A11}
\rho \mho =\bigl(\e\otimes (\rho{\restr}\cal{L}_*)\bigr).
\end{equation}

The statement ({\it i\/}) now follows from the definition
of $\Gamma_{\inv}$ and from the
facts that $\e(\lambda )=\mu +\mu^{-1}$ and $p_{kn}
(\mu +\mu^{-1} )\neq 0$. \end{pf}

Using \eqref{A6} and the definition of $\delta$, it can be shown that
\begin{equation}\label{A12}
\delta (\Gamma^n_{\inv} )\subseteq \sideset{}{^\oplus_*}\sum_{ij}
(\Gamma^i_{\inv}\otimes \Gamma^j_{\inv} ), \end{equation}
for each $n\in\Bbb{N}$, where the sum is taken over pairs $(i,j)$ satisfying
$\vert i-j\vert \leq n\leq i+j$. In particular,
\begin{gather*}
\delta (\vartheta )^{0,n} =d_n \varepsilon \otimes \vartheta \\
\delta (\vartheta )^{n,0} =d_n \vartheta \otimes \varepsilon ,
\end{gather*}
for each $\vartheta \in\Gamma^n_{\inv} ,$ with
$d_n \in\Re\setminus\{0\}$. This implies that singlet
components of $\omega$ do not figure in nonsinglet
components of the curvature.

\end{document}